\documentclass[11pt]{amsart}
\def\Ga{\Gamma}
\def\ga{\gamma}

\def\De{\Delta}
\def\de{\delta}
\def\La{\Lambda}
\def\la{\lambda}
\def\eps{\varepsilon}
\def\ka{\kappa}
\def\om{\omega}
\def\Om{\Omega}

\def\dom{\operatorname{dom}}

\def\bb{\bar{b}}
\def\bx{\bar{x}}

\def\Ker{\operatorname{Ker}\,}
\def\Im{\operatorname{Im}\,}

\def\dimt{\operatorname{dim}_\tau}

\def\sup{\operatorname{sup}}
\def\inf{\operatorname{inf}}

\def\C{\mathbb C}
\def\R{\mathbb R}
\def\Z{\mathbb Z}

\def\E{\mathcal E}

\def\Tr{\operatorname{Tr}}

\def\End{\operatorname{End}}

\def\A{{\mathcal A}}
\def\B{{\mathcal B}}

\def\F{{\mathcal F}}

\def\cl{\operatorname{cl}}
\def\H{\mathcal H}
\def\h{\mathbb H}

\def\E{\mathcal E}
\def\U{{\mathcal U}}
\def\<{\langle}
\def\>{\rangle}

\newtheorem{theorem}{Theorem}[section]
\newtheorem{lemma}[theorem]{Lemma}
\newtheorem{proposition}[theorem]{Proposition}
\newtheorem{corollary}[theorem]{Corollary}

\theoremstyle{definition}
\newtheorem{definition}[theorem]{Definition}

\theoremstyle{remark}
\newtheorem{remarks}[theorem]{Remarks}

\numberwithin{equation}{section}

\begin{document}

\title[Twisted $L^2$ invariants and asymptotic Morse inequalities]
{TWISTED L$^2$ INVARIANTS OF NON-SIMPLY CONNECTED MANIFOLDS AND ASYMPTOTIC
L$^2$ MORSE INEQUALITIES.}
\author{Varghese Mathai}
\address{Department of Mathematics, University of Adelaide, Adelaide 5005,
Australia}
\email{vmathai@maths.adelaide.edu.au}
\author{Mikhail Shubin}
\address{Department of Mathematics, Northeastern University, Boston, Mass., USA}
\email{shubin@neu.edu}

\dedicatory{ Dedicated to Mark Iosifovich Vishik on the occasion of his
75th birthday}
\date{August 20, 1996}
\subjclass{Primary: 58G11, 58G18 and 58G25.}
\keywords{Twisted L$^2$ Betti functions, twisted Novikov-Shubin functions,
Morse theory of closed 1 forms, asymptotic $L^2$ Morse inequalities, von
Neumann algebras.}

\begin{abstract}
 We develop the theory of {\em twisted} $L^2$-cohomology and
twisted spectral invariants for flat  Hilbertian bundles over compact manifolds.
They can be viewed as functions on $H^1(M, {\mathbb R} )$ and they
generalize the
standard  notions. A new feature of the twisted $L^2$-cohomology theory
is that in addition to satisfying the  standard $L^2$ Morse inequalities,
they also satisfy
certain {\em asymptotic} $L^2$ Morse inequalities. These reduce to the
standard Morse inequalities in the finite dimensional case, and when the
Morse 1-form
is exact.
We define the {\em extended} twisted
$L^2$ de Rham cohomology and prove
the {\em asymptotic} $L^2$ Morse-Farber inequalities, which give
quantitative lower bounds for the Morse numbers of a Morse 1-form on $M$.
\end{abstract}

\maketitle

%------------------------------------------------------------
%------------------------------------------------------------
\section*{Introduction}

Let $M$ be a  compact closed $C^\infty$-manifold. A $C^\infty$-function
$f:M\to\R$
is called {\it Morse function} if any critical point $x$ of $f$ (i.e. point
$x\in M$
such that $df(x)=0$) is non-degenerate. This means that the determinant of
the Hessian
does not vanish, i.e.
$$\det\left(\frac{\partial^2f(x)}{\partial x_i\partial x_j}\right)\not=0\,,$$
where the partial derivatives are taken in local coordinates.
It follows that  all critical points are isolated, therefore there
is only finite number of them. The {\it index} of the critical point $x$ is
the number
of the negative eigenvalues of the Hessian. Denote $m_j=m_j(f)$ the number
of critical points
with the index $j$.

M.Morse discovered a connection between the behavior of $f$ and topology of $M$.
In particular, the numbers  $m_j$ (which are  called {\it Morse numbers})
are related
with the (real) Betti numbers $b_j=b_j(M)$ by the {\it Morse inequalities}
$$m_j(f)\ge b_j(M),\quad j=1,\dots,n, \leqno (0.1)$$
where $n=\dim_\R M$. There are also more general inequalities
$$\sum_{p=0}^k (-1)^{k-p}m_p\ge \sum_{p=0}^k (-1)^{k-p}b_p,\quad k=1,\dots
n.\leqno (0.2)$$
 Note that adding two
inequalities (0.2) with $k=j-1$ and $k=j$, we obtain (0.1).

The Morse inequalities can be applied to estimate the Morse numbers $m_j$
from below.
For example, it follows from (0.1) that any Morse function of the
2-dimensional torus
$T^2=\R^2/\Z^2$ should have at least 2 saddle points (critical points of
index 1).
In other direction, knowing explicitly a Morse function on $M$, we can estimate
its Betti numbers from above. This has important applications e.g. in
complex analysis.

We refer to \cite{Mi} for an excellent exposition of the classical Morse
theory and its
applications.

The Morse theory has numerous generalizations, developments and
applications. We will
only discuss  the directions which are most relevant for this paper.

E.Witten \cite{Wi} suggested a new proof of the Morse inequalities, which
is completely
analytic. He suggested to deform the de Rham complex by replacing the de
Rham external
differential by a deformed differential
$$d_s=\exp(-sf)d\exp(sf)=d+se(df), \quad s\gg 0, \leqno (0.3)$$
where $f$ is a Morse function on $M$, $e(df)$ is the operator of the
external multiplication
of forms by the 1-form $df$. Though the dimensions of the cohomology spaces
do not change
under this deformation, the Laplacians acquire a big parameter $s$ and it
becomes possible
to apply semiclassical asymptotics with $h=1/s$. The explicit form of the
deformed Laplacians
shows that their eigenfunctions with small eigenvalues become localized
near the critical points of $f$. The number of small eigenvalues
(multiplicities taken into account) in forms of degree $j$ can be calculated
and it equals $m_j$. Since $0$ is among these eigenvalues, this implies the
inequalities
(0.1), because  the de Rham cohomology space is isomorphic to the space
of harmonic forms of the corresponding degree due to the Hodge theory.

S.Novikov \cite{{Nov},{Nov2}} suggested to replace the Morse function $f$
by a closed
Morse 1-form on $M$. Such a  form can be considered as ``multivalued Morse
function" on $M$.
This theory was further developed by M.Farber \cite{F2} and A.Pazhitnov
(\cite {{P},{P2}}).

It is quite natural to consider more general deformations than (0.3). In
particular,
we can take the deformation
$$d_s=d+se(\om), \quad s\gg 0, \leqno (0.4)$$
where $\om$ is a 1-form which is not necessarily exact. If $\om$ is closed,
we arrive to
a reinterpretation of the Novikov theory (with ``multivalued Morse
functions"), which
was used by M.Farber and A.Pazhitnov. However it is also interesting to
take $\om$ which is
not closed, in spite of the fact that then $d_s^2\not=0$. The case when
$\om$ is dual to a Killing vector field, was considered by Witten
\cite{Wi}. A more
general situation which leads to Morse-type inequalities for arbitrary
vector fields was
studied by Novikov (see Appendix to \cite{NS} and also \cite{Sh3}).

It was noticed by S.Novikov and M.Shubin \cite{NS} that the Betti numbers
$b_j$ in (0.1)
and (0.2) can be replaced by so-called $L^2$ Betti numbers $\bar b_j$ (or
von Neumann Betti
numbers) which were introduced by M.Atiyah \cite{A}.  They can be defined as
von Neumann dimensions (associated with the fundamental group $\pi_1(M)$) of
the spaces of $L^2$ harmonic forms on the universal covering ${\widetilde
M}$ of $M$,
where ${\widetilde M}$ is considered with a Riemannian metric which is
lifted from $M$
(or, equivalently, a Riemannian metric which is invariant under the action of
$\pi_1(M)$ on ${\widetilde M}$ defined by the deck transformations).
J.Dodziuk \cite{Do}
proved that the $L^2$ Betti numbers are homotopy invariants of $M$.

The $L^2$ Morse inequalities can be applied e.g. to prove that
some topological requirements, imposed
on $M$, imply existence of non-trivial $L^2$ harmonic forms on ${\widetilde
M}$ (\cite{NS}).

The Witten method can be applied to the $L^2$ situation as well (\cite{Sh})
in spite
of the fact that the universal covering is generally non-compact and the
Laplacians
have continuous spectrum. The main tool is a variational priciple (see e.g.
\cite{ES})
which works for von Neumann dimensions.

In $L^2$ situation Morse theory is also naturally connected with the
spectrum-near-zero phenomenon in topology which was discovered by
S.Novikov and M.Shubin (see \cite{NS2} and also \cite{{ES}, {E}, {GS},
{Lo}, {LL}, {F}}).
M.Gromov observed that if
$0$ is in the spectrum of the Laplacian on $p$-forms on ${\widetilde M}$, then
$m_p(f)>0$ for any Morse function $f:M\to\R$, in spite of the fact that it
might happen
that in this situation $b_p={\bar b}_p=0$, so
the positivity of $m_p$ does not follow neither from the classical Morse
inequalities,
nor from their $L^2$-version. A quantitative version of this observation
was given by M.Farber \cite{F} with the help of his extended cohomology theory
which puts the spectrum-near-zero invariants into a cohomological context.

A more general context for $L^2$-cohomology is a flat vector bundle $\E$ on
a compact
manifold $M$, such that the fiber $E$ is a Hilbert module over a finite von
Neumann
algebra $\A$. Such a bundle is called {\it Hilbertian bundle}.
The $L^2$-functions and $L^2$-forms on the universal covering of $M$ can be
interpreted
as sections of such a bundle with the fiber $\ell^2(\pi)$, where $\pi=\pi_1(M)$
is considered as a discrete group, $\ell^2(\pi)$ is the Hilbert space of
square-summable complex-valued functions on $\pi$.

In this paper we study the relationship between the
Morse theory of closed 1-forms on $M$ and a {\em twisted}
$L^{2}$-cohomology of a
flat Hilbertian bundle over $M$. Here ``twisted" means that the covariant
derivative
of the flat connection is deformed by adding $e(\om)$, where $\om$
is a closed form on $M$.

A new feature of the twisted $L^{2}$-cohomology  is that in
addition to  the  standard $L^2$ Morse inequalities we also have
{\em asymptotic} $L^2$ Morse inequalities. Here the big parameter is provided
by the Witten-type deformation of the type (0.4)
(with a closed 1-form $\om$) for the covariant derivative.
We also study the twisted
analogue of the spectrum-near-zero invariants. All these invariants can be
viewed as functions on $H^1(M, {\mathbb R})$ possessing some upper
semi-continuity properties.

 We begin with a review of
some background material
on Hilbertian modules over a von Neumann algebra and flat Hilbertian
bundles over a compact manifold $M$ in
section 1. In section 2 we introduce a {\em twisted} analogue
of all  $L^2$-invariants
associated with a flat Hilbertian bundle $\E\to M$ over $M$,
and briefly discuss their
main properties.
%We use these
%twisted $L^2$-invariants to prove the asymptotic $L^2$ Morse
%inequalities in the next section.
In section 3 we review the Morse theory of closed 1-forms on $M$
and  prove the standard-type $L^2$ Morse inequalities (see \cite{NS}
and \cite{Sh} for the case of the regular representation),
as well as the {\em asymptotic} $L^2$ Morse inequalities, and we also give
some applications. The proofs are based on an analogue of the Witten
deformation technique, which has also been applied in the $L^2$-context
by Burghelea, Friedlander, Kappeler,  McDonald
\cite{BFKM} and Shubin \cite{Sh}, in situations which are somewhat different
to those considered here.
In section 4 we briefly review virtual Hilbertian modules as objects in
the extended category introduced by Farber \cite{F}. Then, acting
similarly to \cite{Sh2},
we define the extended twisted $L^2$ de Rham cohomology. The main result here is
the {\em asymptotic} $L^2$ Morse-Farber inequalities, which give
quantitative lower bounds for the Morse numbers of a Morse 1-form on $M$.
In section 5 we end with some calculations and further applications.

On completing a preliminary version of our paper, we received a preprint
of M. Braverman and M. Farber \cite{BF}, which also proved the asymptotic
standard $L^2$ Morse inequalities, but only for the special case of
{\em residually finite} fundamental groups, and using
in an essential way the results of L\"uck \cite{Lu}. In this case, they also
prove the degenerate asymptotic $L^2$ Morse inequalities.

\section{Preliminaries}

In this section we establish the main notation of the paper, and review some
basic facts about Hilbertian $\A$-modules. (See \cite{{F}, {CFM}} for
details.)
We refer to \cite{{Di},{T}} for the
necessary definitions on von Neumann algebras.

Let ${\A}$ be a finite  von Neumann algebra
with a fixed finite,
normal and faithful trace $\tau:{\A}\rightarrow \C$. We will always assume
that this
trace is normalized i.e. $\tau(1)=1$.
The involution
in ${\A}$ will be denoted $*\,$.

By $\ell^{2}({\A})$ we denote the
completion of ${\A}$ with respect to the scalar product
$\langle a,b\rangle=\tau(b^{*}a)$, for $a,b\in{\mathcal A}$.
For example, let $\Gamma$ be a finitely generated discrete group and
$\ell^{2}({\Gamma})$ denote the Hilbert space of square summable
functions on $\Gamma$. Let $\rho$ denote the left regular representation
of $\Gamma$ on $\ell^{2}({\Gamma})$. This extends linearly to a representation
of the group algebra ${\mathbb C}(\Gamma)$. The weak closure of
$\rho({\mathbb C}(\Gamma))$ is called the group von Neumann
algebra, denoted by $ {\mathcal U}(\Gamma)$. The trace $\tau$ is given
by evaluation at the identity element of $\Gamma$, i.e.
$$\tau(a)=\langle a\de_e,\de_e\rangle,\ \ a\in\U(\Ga),$$
where $e$ is the neutral element in $\Ga$, $\de_e\in\ell^2(\Ga)$,
$\de_e(x)=1$ if
$x=e$, and $0$ otherwise.

Recall that a {\em Hilbert module} over ${\A}$
is a Hilbert space $M$
together with a continuous left ${\mathcal A}$-module structure such that there
exists an isometric ${\mathcal A}$-linear embedding of $M$ into
$\ell^{2}({\mathcal A})\otimes H$, for some Hilbert space $H$.
Note that this embedding is not
part of the structure.  A Hilbert module $M$ is {\it finitely generated} if it
admits an embedding $M\rightarrow\ell^{2}({\mathcal A})\otimes H$ as above with
finite-dimensional $H$.
Note that a Hilbert module comes with a particular scalar product.

 A {\em Hilbertian module}
is a topological vector space $M$
with continuous left ${\mathcal A}$-action such that there exists a scalar
product
$\langle\;,\;\rangle$ on $M$ which generates the topology of $M$ and such
that $M$ together
with $\langle\;,\;\rangle$ and with the ${\mathcal A}$-action is a Hilbert
module.

If $M$ is a Hilbertian module,
then any scalar product $\langle\;,\;\rangle$ on
$M$ with the above properties will be called {\em admissible}.
It can be proved that any other
choice of an admissible scalar product gives an isomorphic Hilbert module
(\cite{CFM})
but such choice introduces an additional structure.  The situation here is
similar to the case of finite-dimensional vector spaces: any
choice of a scalar product on a vector space
produces an isomorphic Euclidean vector space.

Let $M$ be a Hilbertian module and let $\langle\;,\;\rangle$ be an admissible
scalar product. Then $\langle\;,\;\rangle$ must be compatible with the
topology on $M$ and with the ${\mathcal A}$-action.
The last condition means that
the involution on ${\mathcal A}$ determined by the scalar product $\langle\;,
\;\rangle$ coincides with the involution of the von Neumann algebra
${\mathcal A}:$
$$
\langle\lambda\cdot v,w\rangle=\langle v,\lambda^{*}\cdot w\rangle
$$
for any $v,w\in M$, $\lambda\in{\mathcal A}$.

If $\<\ ,\ \>$ and $\<\ ,\ \>_1$ are two admissible
scalar products on a Hilbertian module $M$  then the Hilbert modules
$(M, \<\ ,\ \>)$ and $(M, \<\ ,\ \>_1)$ are isomorphic.
Therefore we can define
{\em finitely generated Hilbertian
modules} as those for which the corresponding Hilbert modules
(obtained by a choice of an admissible scalar product) are finitely
generated. Note that the {\em von Neumann dimension}
of a Hilbertian module $M$, denoted $\dim_{\tau}(M)$, is also well defined
(we will recall the definition later).

A {\em morphism} of Hilbertian modules is a continuous linear map
$f:M\to N$, commuting with the $\A$-action.
Note that the kernel of any morphism $f$ is again a Hilbertian module. Also,
the closure of the image $\cl(\Im(f))$ is a Hilbertian module.

Let $\B=B_{\mathcal A}(M)$
denote the set of endomorphisms of $M$ as  Hilbertian module. It can also
be described
as the commutant of the action of $\A$ on $M$, so we will sometimes refer to it
simply as the {\em commutant}.

Any choice of an admissible scalar product $\langle\;,\;\rangle$ on $M$
defines obviously a $*$-operator on $\B$ (by assigning to an operator
its adjoint) and turns $\B$ into a von Neumann algebra.  Note that this
involution $*$ depends on the scalar product $\<\ ,\ \>$ on $M$;
if we  choose another admissible scalar product $\<\ ,\ \>_1$
on $M$ then the new involution will be given by
$$
f\mapsto \ A^{-1}f^*A\qquad \text{for}\quad f\in \B,
$$
where $A\in\B$ is the operator defined by
$\<v,w\>_1\ =\ \<Av,w\>$ for $v,w\in M$.

If $M$ is finitely generated, then
the trace $\tau:{\mathcal A}\rightarrow \C$
determines canonically a trace on the commutant
$$
\Tr_{\tau}:\B\ =\ \B_{\A}(M)\rightarrow \C
$$
which is finite, normal, and faithful.

We now briefly describe this trace.
Suppose first that $M$ is {\em free},
that is, $M$ is isomorphic
to $l^2(\A)\otimes\C^k$ for some $k$. Then the commutant
$\B$ can be identified with the set of all
$k\times k$-matrices over the algebra $\A$, acting from the right on
$l^2(\A)\otimes\C^k$ (the last module is viewed as the set of row-vectors
with components in $l^2(\A)$). If $\alpha\in\B$ is an element represented
by a $k\times k$ matrix $(\alpha_{ij}),$ then one defines
$$ \Tr_{\tau}(\alpha)\ =\ \sum_{i=1}^k \ \tau(\alpha_{ii})$$
This formula gives a trace on $\B$ which is finite, normal and faithful.
Note also that this trace $\Tr_{\tau}$
does not depend on the representation of $M$
as the product $l^2(\A)\otimes\C^k$.

Suppose now that $M$ is an arbitrary finitely generated Hilbertian module.
Then $M$ can be embedded into a free module $F$ as a direct summand.
Let $\pi_M:F\to M$ and $i_M:M\to F$ denote the projection and embedding
correspondingly.
Then for $f\in\B_{\A}(M)$ the formula
$$\Tr_\tau(f)\ =\ \Tr_\tau(i_M\circ f\circ \pi_M) $$
defines a trace
$$\Tr_{\tau}: \B_{\A}(M)\to \C,$$
satisfying
$$\Tr_{\tau}(fg) \ = \ \Tr_{\tau}(gf)$$
for all $f,g\in \B_{\A}(M)$.
This trace is independent of all choices. In particular we can define
$\A$-dimension of $M$
by the formula
$$\dimt M=\Tr_\tau(\hbox{Id}_M)\,.$$

\subsection{Hilbertian $\A$-complexes}
 Let us consider a sequence
$$L_\bullet:\;0\longrightarrow L_0\buildrel{d_0}\over\longrightarrow
L_1\buildrel{d_1}\over\longrightarrow \ldots
\buildrel{d_{n-1}}\over\longrightarrow L_n\longrightarrow 0,\leqno (1.1)$$
where all $L_j,\ j=1,\dots,n,$ are Hilbertian $\A$-modules,
$d_j$  are closed densely defined linear operators which commute with the
action of $\A$ in the sense that $d_ja=ad_j$ on $\dom(d_j)$ for all $a\in\A$.
(Here $\dom(d_j)$ denotes the domain of $d_j$ and it is a dense linear
subspace in $L_j$.)
In particular this means that $a(\dom(d_j))\subset\dom(d_j)$ i.e. the
domains $\dom(d_j)$
are $\A$-invariant.

Such a sequence  is called a {\it complex of Hilbertian} $\A$-{\it modules} if
 $d_jd_{j-1}=0$ on $\dom(d_{j-1})$ for all $j$. In particular this means that
$\Im d_{j-1}\subset \Ker d_j$. Note that $\Ker d_j$ is always closed.

We will call a complex (1.1) {\em finite} if all Hilbertian modules
$L_j$ are finitely generated and all differentials  $d_j$ are (bounded)
morphisms
of Hilbertian $\A$-modules.

The reduced $L^2$-cohomology groups of $L_\bullet$ are defined as
the Hilbertian $\A$-modules
$$H^p_{(2)}(L)=\frac{\Ker d_p}{\cl(\Im d_{p-1})},\ \ p=0,1,\dots,n. $$
(Here by definition $d_{-1}$ and $d_n$ are zero morphisms.) Denote
$m_p=\dimt L_p$, $\bb_p=\dimt H^p_{(2)}(L).$

The following Lemmas are well known:

\begin{lemma} Suppose that $m_j<\infty$ {\it for all} $j$. {\it Then}
$$\sum_{j=0}^n (-1)^j m_j=\sum_{j=0}^n (-1)^j\bb_j.  $$
\end{lemma}

\begin{lemma} {\it If} $m_j<\infty$ {\it for all} $j$ {\it then}
$$\sum_{j=0}^p (-1)^{p-j} m_j\geq\sum_{j=0}^p (-1)^{p-j}\bb_j  $$
{\it for every} $p=0,\dots n$.\end{lemma}

The proofs do not differ from the proofs of the corresponding
statements for the von Neumann algebra $\A={\mathbb C}$ (and for the spaces
$L_j$ which are
finite-dimensional in the usual sense) except almost isomorphisms should be
used instead of usual isomorphisms (see \cite{Sh} for more details).

If $M_\bullet$ and ${N}_\bullet$ are Hilbertian $\A$-complexes, then a
{\em morphism of Hilbertian $\A$-complexes}
$f: {M}_\bullet \rightarrow {N}_\bullet$ is a sequence
$f_{k} : {M}_{k}\rightarrow N_{k}$ of morphisms of Hilbertian $\A$-modules
such that $\displaystyle   f_{k+1} d_{k}w = d_{k} f_{k} w$
for all $w \in\dom (d_{k})$.  A {\em homotopy} between two morphisms
$\displaystyle f, g: {M}_\bullet \rightarrow {N}_\bullet$
is a sequence of morphisms of Hilbertian $\A$-modules
${T}_{k} :  M_{k} \rightarrow N_{k-1}$ such that
$\displaystyle f_{k} - g_{k} = {T}_{k+1} d_{k} + d_{k-1}
{T}_{k}$ on $\dom(d_{k})$.  Homotopy is an equivalence relation.

Two Hilbertian $\A$-complexes ${M}_\bullet$ and ${N}_\bullet$ are said to be
homotopy
equivalent if there exist morphisms $f:{M}_\bullet \rightarrow {N}_\bullet$ and
$g: {N}_\bullet \rightarrow {M}_\bullet$ such that $fg$ and $gf$
are homotopic to the
identity morphisms of ${N}_\bullet$ and ${M}_\bullet$ respectively.  Homotopy
equivalence is an equivalence relation.

If ${M}_\bullet$ is a Hilbertian $\A$-complex,
then we can define
functions $F_{k} (\lambda, {M})$ as
$$
F_{k} (\lambda, {M}) \equiv
\sup\left( \dim_{\tau} L : L \in  {\mathcal S}^{(k)}_{\lambda}
({M}) \right)
$$
where ${\mathcal S}^{(k)}_{\lambda} ({M})$ denotes the set of all closed
$\A$-invariant subspaces of ${M}_{k} / \ker d_{k}$ such that
$ L \subset\hbox{dom}(d_{k}) /  \ker d_{k}$ and
$\displaystyle \| d_{k} w \| \leq
\sqrt{\lambda} \| w \| $ for $w \in L$.
(The norm on the right hand side is the
quotient norm). Then $\la\mapsto F_{k} (\la,M) $ is an
increasing function on ${\mathbb R} $ with
values in $[0, \infty]$ and $F_{k}(\lambda, M) = 0$ if
$\lambda < 0$.

Given a Hilbertian $\A$-complex $M_\bullet$,  consider
the {\em Laplacian} $\Delta_{k} \equiv d_{k-1} {\delta}_{k-1} +
{\delta}_{k} d_{k}$, where ${\delta}_{k}$ denotes
the $L^{2}$ adjoint of $d_{k}$. It is a self-adjoint operator in $M_k$
and it has the spectral decomposition
$$
\displaystyle \Delta_{k} = \int_{0}^{\infty} \lambda  d E_{\lambda}.
$$
Then the von Neumann spectral density function is defined as

$$
\displaystyle {N}_{k}(\lambda, {M}) =\Tr_{\tau}E_\la,
%=\int_{0-}^{\lambda}\Tr_{\tau} (d E_{\mu}),
$$
and can be expressed through the functions $F_k$ as follows (\cite{GS})
$$
\displaystyle {N}_{k}(\lambda, {M}) =
F_{k-1}(\lambda, {M})+F_{k}(\lambda, {M})
+ b_{(2)}^{k}({M}),
$$
where $b_{(2)}^{k}({M})= \dim_{\tau} \Ker \Delta_{k}$ are the $L^2$
Betti numbers.

Two functions $F(\lambda)$ and $G(\lambda)$ on $(0,\infty)$ satisfy
$F \prec\!\!\prec G$ if there exist  positive constants
$C , \la_0$ such that
$F(\lambda) \leq G(C \lambda)$ for all $\la\in(0,\la_0)$.

If $F \prec\!\!\prec G$ and $G \prec\!\!\prec F$, then
$F$ and $G$ are said to
be dilatationally equivalent and we write $F \sim G$.
Roughly speaking, in this case the small
$\lambda$    asymptotics of $F(\lambda)$ and $G(\lambda)$ are the same.
The following basic abstract theorem is due to Gromov and Shubin
\cite{GS}:

\begin{theorem}
 Let $f:{M} \rightarrow {N}$ and $g: {N} \rightarrow
{M}$ be morphisms of Hilbertian $\A$-complexes such that $gf$ is
homotopic
to the identity morphism of ${M}$.  Then $F_{k}(\lambda, {M})
\prec\!\!\prec
F_{k}(\lambda, {N})$ and
$b_{(2)}^{k}({M}) \leq b_{(2)}^{k}({N})$
for all $k$.  Hence if ${M}$ and ${N}$ are
homotopy equivalent, then $F_{k}(\lambda, {M})
\sim F_{k}(\lambda, {N})$ and $b_{(2)}^{k}({M}) = b_{(2)}^{k}({N})$.

\end{theorem}

\subsection{Flat Hilbertian $\A$-Bundles} Let $E$ be a finitely generated
Hilbertian $(\A-\pi)$-bimodule. This means first that $\A$ acts on $E$
from the left,  so that with respect to this action $E$ is a finitely generated
Hilbertian module; and second that  $\pi$ is a discrete group acting on $E$
from the right and the action of $\pi$ commutes with that of $\A$.
A Hilbertian $(\A-\pi)$-bimodule $E$ is said to be {\em unitary} if there exists
an admissible scalar product $\<\ ,\ \>$ on $E$ such that the action of
$\pi$ on $E$ preserves this scalar product.

Let $M$ be a connected, closed, smooth
manifold with fundamental group $\pi$. Let $\widetilde M$ denote the
universal covering of $M$.
{\it A flat Hilbertian $\A$-bundle with fiber $E$
over $M$} is an associated bundle $p:{\E}\to M$. This means that
$\E\ =\ (E\times\widetilde M)/\sim$ with its natural projection onto $M$, where
 $(v,x) \sim (vg^{-1},gx)$ for all $g\in \pi$, $x\in \widetilde M$ and
$v\in E$.
Then $p:{\E}\to M$ is a locally trivial bundle of topological vector spaces
(cf. chapter 3 of \cite{La}) which has a natural fiberwise left action of $\A$.

A flat Hilbertian $\A$-bundle $\E\to M$ over $M$, with fiber $E$, is said
to be {\em unitary} if the Hilbertian $\A$ -module $E$ is unitary.

Any smooth section $s$ of $\E\to M$ can be uniquely represented by
a smooth equivariant map $\phi:\widetilde M\to E$, where ``equivariant"
means that
$\phi(gx)=\phi(x)g^{-1}$ for all $g\in\pi$ and $x\in \widetilde M$. Given
such a map
$\phi$, the corresponding section assigns to every $y\in M$, the
equivalence class of the pair $(x,\phi(x))$, where $x$ is a lifting
of the point $y$.

Given a flat Hilbertian
$\A$-bundle $\E\to M$ over a closed connected manifold $M$,
one can consider the space of smooth differential $j$-forms on $M$
with values in $\E$; this space will be denoted by $\Omega^j(M,\E)$.
It is naturally
defined as a left $\A$-module and can be written as
$$\Omega^j(M,\E)=C^\infty(M, \E)\bigotimes_{C^\infty(M)}
\Omega^j(M)\;,\leqno (1.2)$$
where $C^\infty(M,\E)$ is the set of all $C^\infty$-sections of $\E$ over $M$.
An element of $\Omega^j(M,\E)$
can be also uniquely represented as a $\pi$-invariant element in
$$C^\infty(\widetilde M, E)\bigotimes_{C^\infty({\tilde M})}
\Omega^j(\widetilde M)\leqno (1.3)$$
with respect to the total (diagonal) action of $\pi$.
Here $C^\infty({\widetilde M},E)$ is the space of $E$-valued
$C^\infty$-functions on
$\widetilde M$, and $\Om^j({\widetilde M})$ is the space of
$C^\infty$-forms of degree $j$
on ${\widetilde M}$.

% and $\pi$-invariance is taken with respect to  the total (diagonal)
%$\pi$ action,
%that is, the tensor product of the actions of $\pi$ on $\Omega^j(\widetilde M)$
%and on   $C^\infty(\widetilde M, E)$. More precisely, if
%$\omega\in\Omega^j(\widetilde M)$ and
%$f\in C^\infty({\widetilde M}, E)$,
%then
%$f\otimes\omega$ is said to be $\pi$-invariant if
%$f(gx)g\otimes (g^\ast\omega)(x)= f(x)\otimes\omega(x)$
%for all $x\in\widetilde M$ and all $g\in\pi$.

 An  {\it $\A$-linear connection}
on a flat Hilbertian $\A$-bundle $\E$ is defined as an $\A$-homomorphism
$$\nabla : \Omega^j(M,\E)\to \Omega^{j+1}(M,\E)$$
which is given for all $j$ and satisfies the Leibniz rule
$$\nabla(f\omega ) = df\wedge\omega + f\nabla (\omega)$$
for any $\A$-valued function $f$ on $M$ and
for any $\omega\in \Omega^j(M,\E)$. This connection is called {\it flat} if
$\nabla^2=0$.
On a flat Hilbertian $\A$-bundle $\E$, as defined above,
there is a {\it canonical flat
$\A$-linear connection} $\nabla$ which is given as follows: under the
identification of $\Omega^j(M,\E)$ given in the previous paragraph,
one defines the connection $\nabla$ to be
the de Rham exterior derivative acting on the second factor in (1.3).

\subsection{Hermitian metrics and L$^2$ scalar products}
{\it A Hermitian metric} on a flat Hilbertian $\A$-bundle $p:\E\to M$
is a smooth family of admissible scalar products on the fibers. Any
Hermitian metric on $p:\E\to M$ defines a wedge-type product
$$
\wedge : \Omega^p(M,\E)\otimes \Omega^{q}(M,\E)\rightarrow
\Omega^{p+q}(M)
$$
similar to the finite dimensional case (see e.g. (1.20), Chapter 3, in
\cite{We}).
Note that this product is antilinear with respect to the second factor, and
it is
Hermitian-antisymmetric i.e.
$$\alpha\wedge\beta=(-1)^{pq}\overline{\beta\wedge\alpha},\quad
\alpha\in\Om^p(M,\E),\ \beta\in\Om^q(M,\E)\;.$$

Suppose we are also given a Riemannian metric on $M$.
Then we can define the Hodge star-operator
$$\ast:\Omega^j(M,\E)\to \Omega^{n-j}(M,\E\otimes o(M)),$$
where $o(M)$ is the orientation bundle of $M$ which is a real line-bundle
(trivial if
$M$ is oriented). The operator $\ast$ is
a complex linear operator defined as the complexification of the
standard real Hodge star operator. It acts on form-coefficients (without
affecting the fiber)
i.e.
$$\ast(f\otimes\om)=f\otimes(\ast\om),\quad f\in C^\infty(M,\E),\
\om\in\Om^\bullet(M)\,. $$

The Hermitian metric on $p:\E\to M$ together with a Riemannian metric
on $M$ determines a scalar product on $\Omega^i(M,\E)$ in the
standard way. Namely, one sets
$$(\omega,\omega^\prime)\ =\ \int_M \omega\wedge\ast\omega^\prime $$
(cf. \cite{We}, Section 2 in Chapter 5).
%Generally one can use partition of unity subordinated to the covering
consisting
%of oriented pieces.

With this scalar product $\Omega^i(M,\E)$
becomes a pre-Hilbert space. Define the space of
$L^2$ differential $j$-forms on $M$ with coefficients in
$\E$, which is denoted by $\Omega_{(2)}^j(M,\E)$, to be the Hilbert
space completion of $\Omega^j(M,\E)$. We will tend to ignore the
scalar product on $\Omega_{(2)}^j(M,\E)$ and view it as a
Hilbertian $\A$ module (not necessarily finitely generated).

The connection $\nabla$ on $\E$ extends by closure
to a closed, unbounded, densely defined operator
$\nabla:\Omega_{(2)}^j(M,\E)\to\Omega_{(2)}^{j+1}(M,\E)$.

\subsection{Reduced L$^2$ cohomology}
Given a flat Hilbertian $\A$ bundle $p:\E\to M$, one
defines the {\it reduced $L^2$ cohomology with coefficients
in $\E$} as the quotient
$$
H^{j}(M,\E)=\frac{\Ker \nabla/\Omega^{j}_{(2)}(M,\E)}{\cl(\Im\;\nabla/
\Omega^{j-1}(M,\E))}.
$$
Then $H^j(M,\E)$ is naturally defined as a Hilbertian module over $\A$.
The arguments given in \cite{{Do}, {Sh2}}  show that it coincides with the
reduced combinatorial cohomology of
$M$ with coefficients in a locally constant sheaf, determined by $\E$.

\subsection{Hodge decomposition}
The Laplacian $\Delta_j$ acting
on $\E$-valued $L^2$-forms of degree $j$ on $M$ is defined to be
$$
\Delta_j=\nabla \nabla^{*} + \nabla^* \nabla:\Omega_{(2)}^j(M,\E)
\rightarrow\Omega_{(2)}^j(M,\E)\,,
$$
where $\nabla^{*}$ denotes the Hilbert adjoint of $\nabla$ with respect to the
$L^{2}$ scalar product on $\Omega_{(2)}^\bullet(M,\E)$.
It is easy to see that the
Laplacian is a  self-adjoint operator. In fact it can be obtained as the
closure from
the operator given by the same expression on smooth forms.

Let ${\H}^{j}(M,\E)$ denote the closed subspace of  harmonic
$L^2$-forms of degree $j$ with coefficients in $\E$, that is, the kernel of
$\Delta_j$. Note that ${\H}^j(M,\E)$ is a Hilbertian ${\A}$-module.
By elliptic regularity (cf. section 2, \cite{BFKM}), one sees that
${\H}^j(M,\E) \subset \Omega^j(M,\E)$, that is, every $L^2$ harmonic
$j$-form with coefficients in $\E$ is smooth.
Standard arguments then show that one has the following Hodge decomposition
(cf. \cite{Do}; section 4, \cite{BFKM} and also section 3, \cite{GS})
$$
\Omega_{(2)}^j(M,\E) = {\H}^j(M,\E) \oplus \cl(\Im\;\nabla/
\Omega^{j-1}(M,\E)) \oplus \cl(\Im\;\nabla^*/
\Omega^{j+1}(M,\E)).
$$
Therefore it follows that the natural map
$$
{\H}^{j}(M,\E)\rightarrow H^{j}(M,\E)
$$
is an isomorphism Hilbertian ${\mathcal A}$-modules.
The corresponding $L^2$ Betti numbers are
$$
b^{j}_{(2)}(M,\E) = \dimt \left( H^{j}(M,\E)\right).
$$

\begin{definition}
Let $\Delta_j = \int_0^\infty \lambda dE_j(\lambda)$ denote the spectral
decomposition of the Laplacian. The {\it spectral density function} is
defined to be
$N_j(\lambda) = \Tr_\tau(E_j(\lambda))$ and the {\it theta function}
is defined to be
$\Theta_j(t) = \int_{0+}^\infty e^{-t\lambda} dN_j(\lambda)=
\Tr_\tau(e^{-t\Delta_j}) - b^{j}_{(2)}(M,\E)$.

These quantities are well defined because the projection
$E_j(\lambda)$ and the heat operator $e^{-t\Delta_j}$ have smooth
Schwartz kernels which are smooth sections of a bundle over $M\times M$ with
fiber the commutant of $E$, cf. \cite{{BFKM}, {Luk}}.
The symbol $\Tr_\tau$ denotes
application of the canonical trace
on the commutant to the restriction of the kernels to the
diagonal followed by integration over the manifold $M$. This is a trace.

The {\em Novikov-Shubin invariants} are defined as follows (cf. \cite{{ES},
{Lo}, {LL}}):
$$
\alpha_j(M,\E) = \sup\{\beta\in I\!\!R:\Theta_{j}(t)\;\;\mbox{ is }
\;\;O(t^{-\beta})\;\mbox{ as }\;t\rightarrow\infty\}\in[0,\infty]
$$
$$
\overline\alpha_{j}(M,\E)=\inf\{\beta\in
I\!\!R: t^{-\beta}\;\;\mbox{ is }
\;\;O(\Theta_{j}(t)))\;\mbox{ as }\;t\rightarrow\infty\}\in[0,\infty].
$$
\end{definition}

\section{Twisted $L^{2}$ cohomology and twisted $L^{2}$ invariants.}

In this section we introduce the {\em twisted} $L^2$ invariants
associated with a flat Hilbertian $\A$-bundle $p:\E\to M$ over $M$. There
are twisted
analogues of all the invariants which were discussed in the previous
section. They can be viewed as
functions on $H^1(M,{\mathbb R})$ and will be
 used later in the asymptotic $L^2$ Morse
inequalities in the next section.

We will assume that $M$ is a compact Riemannian manifold.
Let $\theta$ be a closed 1-form on $M$.
Consider  the twisted complex
$$
(\Omega_{(2)}^\bullet (M,\E),\nabla_{\theta})
$$
where
$\Omega_{(2)}^\bullet(M,\E)$ denotes the space of $L^{2}$ differential
forms on ${M}$ with coefficients in $\E$, and the differential is given by
$$
\nabla_{\theta}=\nabla+e(\theta)
$$
where $e(\theta)$ denotes exterior multiplication by the 1-form
$\theta$.  Clearly $\nabla^{2}_{\theta}=0$. The closed 1-form
defines a (real) representation of the fundamental group,
$$
\begin{array}{l}
\displaystyle\rho_\theta : \pi_1(M) \to (0,\infty)\\
\displaystyle\rho_\theta (\alpha) = \exp\left( -\int_\alpha \theta\right).
\end{array}
$$
This in turn defines a new flat Hilbertian $\A$-bundle
$p: \E_\theta \to M$ over $M$,
which can be described in several equivalent ways. One way is to note that
the representation $\rho_\theta$ defines a flat real line bundle $p':L_\theta
\to M$ over $M$, where
$L_\theta =\ ({\mathbb R}\times\widetilde M)/\sim$ with its natural projection
onto $M$.
Here $(v,x) \sim (v\rho_\theta (g^{-1}),gx)$ for all $g\in \pi$, $x\in
\widetilde M$ and $v\in
{\mathbb R}$, $\widetilde M$ denotes the universal covering of $M$.
Then the flat Hilbertian $\A$-bundle $\E_\theta$ over $M$, is defined to be
the tensor product $\E \otimes L_\theta$.

The {\em twisted $L^{2}$ cohomology} is defined as
$$
H^{j}_{(2)}(M,\E_\theta)=\frac{\Ker \nabla_
\theta|_{\Omega^{j}_{(2)}(M,\E)}}
{\cl(\Im\nabla_{\theta}|_{\Omega^{j-1}_{(2)}(M,\E)})}\;.
$$
Since $\A$
commutes with the differential
$\nabla_{\theta}$, it follows that $H^{j}_{(2)}(M,\E_\theta)$
is a Hilbertian $\A$-module.
Thus we can define the {\em twisted $L^{2}$ Betti numbers} as
$$
b^{j}_{(2)}(M,\E_\theta)\equiv
\dimt(H^{j}_{(2)}(M,\E_\theta)).
$$
We define the {\em twisted Laplacian}
$$
\Delta_{\theta,j}=\nabla^{*}_{\theta}\nabla_{\theta}+\nabla_{\theta}\nabla^{
*}_{\theta}
$$
acting on $\Omega^{j}_{(2)}(M,\E_\theta)$. Here
$\nabla^{*}_{\theta} = \nabla^* + i(V)$
denotes the formal adjoint of the operator $\nabla_\theta$ and $i(V)$ denotes
contraction with the vector field $V$ which is the Riemannian dual
to the 1-form
$\theta$. Then $\Delta_{\theta,j}$ is a formally
self-adjoint operator with a unique self-adjoint extension, which we denote
by the same symbol.

We denote by ${\mathcal H}^{j}(M,\E_\theta)$ the kernel of
$\Delta_{\theta,j}$,
and we refer to elements in ${\mathcal H}^{j}(M,\E_\theta)$ as {\em twisted
$L^{2}$ harmonic $j$-forms}.
>From the Hodge theorem discussed earlier we know that the Hilbertian
$\A$-modules
${\mathcal H}^{j}(M,\E_\theta)$ and
 $H^{j}_{(2)}(M,\E_\theta)$ are isomorphic.

Let
$$
\Delta_{\theta,j}=\int\lambda\;dE^{j}_{\lambda}(\theta)
$$
denote the spectral decomposition of the twisted Laplacian and
$e^{\theta}_{j}(\lambda,x,y)$ denote the Schwartz kernel of the projection
$E^{j}_{\lambda}
(\theta)$.  Then by elliptic regularity theory,
$e^{\theta}_{j}(\lambda,x,y)$
is smooth in $x,y$ and we can define the twisted spectral distribution function
$$
N_{j}(\lambda,\theta)=Tr_{\tau}(E^{j}_{\lambda}(\theta))=
\int_{M}Tr_{\tau}(e^{\theta}_{j}(\lambda,x,x))dx
$$
where the first $Tr_{\tau}$ denotes the von Neumann trace
on the algebra of $\A$-invariant operators on
$\Omega^{*}_{(2)}({M},\E)$, and the trace under the integral sign means the
trace
in the corresponding fiber. It follows easily that
$$
\lim_{\lambda\rightarrow 0^{+}}N_{j}(\lambda,\theta)=
b^{j}_{(2)}({M},\E_\theta)\;.
$$

\begin{proposition}
\begin{itemize}

\item[(a)] The dilatation class of
$N_{j}(\lambda,\theta)$ depends
only on the cohomology class $[\theta]\in H^1(M,{\mathbb R} )$ of the form
$\theta$.

\item[(b)] $b^{j}_{(2)}(M,\E_\theta)$ depends only on the
cohomology class $[\theta]\in H^1(M,{\mathbb R} )$ of the form $\theta$.
\end{itemize}
\end{proposition}

\begin{proof}   Let $\theta'$ be a closed 1-form on $M$ which is
cohomologous to $\theta$ i.e. $\theta'-\theta=dh$ where $h\in C^\infty(M)$.
 Then since
$$
\nabla_{\theta'}=e^{-{h}}\nabla_{\theta}e^{{h}}\;,
$$
we see that the following diagram of complexes commutes:
$$
\begin{array}{lcccl}
\ldots\rightarrow & \Omega^{j}_{(2)}(M,\E_{\theta'}) &
\stackrel{\nabla_{\theta'}}{\longrightarrow} &
\Omega^{j+1}_{(2)}(M,\E_{\theta'}) & \longrightarrow\ldots \\
 &&&& \\
 & \downarrow\;e^{{h}} & & \downarrow\;e^{{h}} & \\
 &&&& \\
\ldots\rightarrow & \Omega^{j}_{(2)}(M,\E_{\theta}) &
\stackrel{\nabla_{\theta}}{\longrightarrow}
& \Omega^{j+1}_{(2)}(M,\E_{\theta}) & \rightarrow\ldots
\end{array}
$$
and $e^{{h}}$ is a bounded $\A$-invariant map, that is,
$e^{{h}}$
is a morphism of Hilbert $\A$-modules.
Thus by the abstract Theorem 1.3, $N_{j}(\lambda,
\theta)$ and $N_{j}(\lambda,\theta')$ are dilatationally equivalent,
completing the proof of part a).

The proof of part b) follows immediately from part a), by evaluating
the spectral density functions at $s=0$.
\end{proof}

\begin{proposition} $\;$(Poincar\'{e} duality). Let $\E\to M$ be a flat unitary
Hilbertian $\A$-bundle over $M$, and $\theta$ a closed 1-form on $M$. Then
$$
b^{j}_{(2)}({M},\E_{\theta})=b^{n-j}_{(2)}({M},\E_{-\theta}\otimes o(M))\,,
$$
where $n=\dim_{\R}M$, $o(M)$ is the orientation bundle of $M$.

In particular, if $M$ is orientable then
$$
b^{j}_{(2)}({M},\E_{\theta})=b^{n-j}_{(2)}({M},\E_{-\theta})\,.
$$

\end{proposition}

\begin{proof}  Let $\ast$ denote the Hodge star operator, which
induces a linear
isomorphism
$$
\ast:\Omega^{j}_{(2)}({M},\E_{\theta})\rightarrow
C^\infty(M,\E_{-\theta}\otimes \La^{n-j}T^\ast(M)\otimes o(M)).
$$
It is easy to see that $\ast$ intertwines the Laplacians on the bundles
$\E_{\theta}\otimes \La^jT^\ast M$ and
\newline
$\E_{-\theta}\otimes \La^{n-j}T^\ast M\otimes o(M)$.
The result follows.
\end{proof}

We next define the closed subspaces
$$
E_{j}(M,\E_\theta)=\overline{\Im\nabla_{\theta}}\;
\mbox{ in }\;
\Omega^{j}_{(2)}(M,\E_{\theta})
$$
and
$$
E^{*}_{j}(M,\E_\theta)=\overline{\Im
\nabla^{*}_{\theta}}\;\mbox{ in }
\;\Omega^{j}_{(2)}(M,\E_{\theta})\;.
$$
The Kodaira-Hodge decomposition theorem yields
$$
\Omega^{j}_{(2)}(M,\E_{\theta})=E_{j}(M,\E_\theta)
\oplus{\mathcal H}^{j}
(M,\E_\theta)\oplus E^{*}_{j}(M,\E_\theta).
$$
This is an orthogonal decomposition with respect to the scalar product in
$\Omega^{j}_{(2)}(M,\E_{\theta})$.  Since $\Ker \nabla_{\theta}$
is the orthogonal
complement of $E^{*}_{j}(M,\E_\theta)$, we have
$$
\Ker \nabla_{\theta}=E_{j}(M,\E_\theta)\oplus{\mathcal H}^{j}
(M,\E_\theta)\,.
$$
Since $\Ker \nabla^{*}_{\theta}$ is the orthogonal complement of
$E_{j}(M,\E_\theta)$, we have
$$
\Ker \nabla^{*}_{\theta}={\mathcal H}^{j}(M,\E_\theta)\oplus E^{*}_{j}(
M,\E_\theta)\,.
$$
Let $F_j(\lambda,\theta)$ denote the spectral density function of the
operator $\nabla_{\theta}^{*}\nabla_{\theta}$ acting on the subspace $E^{*}_{j}
(M,\E_\theta)$ and $G_{j}(\lambda,\theta)$ denote the spectral
density function of the operator $\nabla_{\theta}\nabla^{*}_{\theta}$ acting on
the subspace $E_{j}(M,\E_\theta)$.
Then clearly one has
$$
N_{j}(\lambda,\theta)=G_{j}(\lambda,\theta)+b^{j}_{(2)}
(M,\E_\theta)+F_{j}(\lambda,\theta)\;.
$$

\begin{lemma}  $F_{j}(\lambda,\theta)=G_{j+1}(\lambda,\theta)$ for
$j=-1,0,1,\ldots,n$, where $n=\dim_{\R} M$ and by fiat
$F_{-1}(\lambda,\theta)=G_{n+1}(\lambda,\theta)=0$.
\end{lemma}

\begin{proof}  Note that
$$
\nabla_{\theta,j}:E^{*}_{j}(M,\E_\theta)\longrightarrow
E_{j+1}(M,\E_\theta)
$$
is an (unbounded) almost isomorphism of Hilbert $\A$-modules,
which intertwines
the operators $\nabla^{*}_{\theta}\nabla_{\theta}$ acting on
$E^{*}_{j}(M,\E_\theta)$ and $\nabla_{\theta}\nabla^{*}_{\theta}$ acting on
$E_{j+1}(M,\E_\theta)$.

Let $U_{\theta,j}$ denote the unitary factor of $\nabla_{\theta,j}$ in its
polar decomposition.  Then
$$
U_{\theta,j}:E_{j}^{*}(M,\E_\theta)\longrightarrow
E_{j+1}(M,\E_\theta)
$$
is a bounded isomorphism of Hilbert $\A$-modules which intertwines the
same pair of operators as before.  This proves the lemma.
\end{proof}

We conclude that
$$
N_{j}(\lambda,\theta)=F_{j-1}(\lambda,\theta)+b^{j}_{(2)}(M,\E_\theta)
+F_{j}(\lambda,\theta).
$$

\begin{lemma}.  We can express
$$
\begin{array}{lcl}
F_{j}(\lambda,\theta)=\sup \Big\{&
\dimt L\;  |\;  L
\mbox{ is a closed subspace of}\;\;\Omega^{j}_{(2)}(M,\E_\theta)/\Ker
\nabla_{\theta} \mbox{ satisfying} \\
&\|\nabla_{\theta}\omega\|\leq \lambda\|\omega\|_{q}\;\;
\mbox{ for all }\;\;\omega\in L \Big\}.
\end{array}
$$
Here the norm $\|.\|$ denotes the usual norm on
$\Omega^{j+1}_{2}(M,\E_\theta)$, whereas $\|.\|_q$ denotes the norm in the
quotient space $\Omega^{j}_{(2)}(M,\E_\theta)/\Ker \nabla_{\theta}\;$.
\end{lemma}

\begin{proof}   This is an immediate consequence of the isomorphism
$$
\frac{\Omega^{j}_{(2)}(M,\E_\theta)}{\Ker \nabla_{\theta}}\cong
E^{*}_{j}(M,\E_\theta)
$$
and the variational principle of \cite{ES}.
\end{proof}

\begin{corollary}.  The dilatation class of $F_{j}(\lambda,\theta)$ is
independent of the choice of Riemannian metric on $M$.\end{corollary}

\begin{proof}   This is immediate from the
expression for $F_{j}(\lambda,\theta)$ obtained in the previous lemma.
\end{proof}

It is also clear that
$$
b^{j}_{(2)}(M,\E_\theta)=\dimt
\left(H^{j}_{(2)}(M,\E_\theta)\right)=\dimt\left(\frac{\Ker \nabla_
\theta|_{\Omega^{j}_{(2)}(M,\E_\theta)}} {\cl(\Im\nabla_{\theta}
|_{\Omega^{j-1}_{(2)}(M,\E_\theta)})}\right)
$$
is independent of the choice of metric on $M$.

\begin{corollary}.  The dilatation class of $N_{j}(\lambda,\theta)$ is
independent of the choice of metric on $M$.\end{corollary}

We define the twisted Novikov-Shubin invariants for $j=0,1,\ldots,\dim M$,
as follows:
$$
\begin{array}{l}
\alpha_{j}(M,\E_\theta)=\sup\{\beta\in
{\mathbb R} :N_{j}(\lambda,
\theta)-b^{j}_{(2)}(M,\E_\theta)\;\mbox{ is}\;O(\lambda^{\beta})
\;\mbox{ as }\;\lambda\rightarrow 0\} \\[+5pt]
{\bar \alpha}_{j}(M,\E_\theta)=
\inf\{\beta\in
{\mathbb R} :\lambda^{\beta}\;\mbox{ is
}\;O(N_{j}(\lambda,\theta)-b^{j}_{(2)}(M,\E_\theta))\;\mbox{ as }\;
\lambda\rightarrow 0\}.
\end{array}
$$
As an immediate consequence of Proposition 2.1 and Corollary 2.6, one has

\begin{corollary}.  $\alpha_{j}(M,\E_\theta)$ and
${\bar\alpha}_{j}(M,\E_\theta)$
are independent of the choice of Riemannian metric on $M$.
Moreover, they depend only on the cohomology class
$[\theta]\in H^1(M,{\mathbb R} )$
of the form $\theta$.
\end{corollary}

Define the {\em twisted von Neumann theta functions} as
$$ \Theta_{j}(t, \theta)= \int_0^\infty e^{-t\lambda} dN_j(\lambda,\theta)
=\Tr_{\tau}(e^{-t\Delta_{\theta,j}}) - b^j_{(2)}(M,\E_\theta).
$$
By the Tauberian theorem relating the large time $t$ asymptotics of the
von Neumann theta function $ \Theta_{j}(t,\theta)$ to the small $\lambda$
asymptotics of the twisted spectral density function
$N_{j}(\lambda,\theta)$ (see Appendix in \cite{GS})
one sees that
$$
\begin{array}{l}
\alpha_{j}(M,\E_\theta)=
\sup\{\beta\in{\mathbb R} :\Theta_{j}(t,\theta)\;
\mbox{ is }\;O(t^{-\beta})\;\mbox{ as }\;t\rightarrow \infty\}, \\[+8pt]
{\bar \alpha}_{j}(M,\E_\theta)=
\inf\{\beta\in{\mathbb R} :t^{-\beta}\;
\mbox{ is}\;O(\Theta_{j}(t,\theta))\;\mbox{ as }\;
t\rightarrow \infty\}.
\end{array}
$$

\section{Asymptotic $L^2$ Morse inequalities.}
In this section, we prove Morse inequalities of two types
for the twisted $L^{2}$ Betti numbers. Except of the standard type inequalities
we also establish
new asymptotic Morse inequalities. We refer the reader to Novikov \cite{Nov2}
and Pazhitnov \cite{P} for discussions on the Morse theory of closed 1-forms
as applied to finite-dimensional flat bundles over compact manifolds.

\subsection{Morse theory of closed 1-forms}
We now briefly review the Morse theory of closed 1-forms on $M$,
which is one way to generalize the Morse theory of $M$
by taking into account
the fundamental group $\pi_{1}(M)$ of $M$.
>From now on we assume that $M$ be a smooth, closed, connected, orientable
manifold and let $p:\widetilde{M}\rightarrow M$ be the universal
cover of $M$.

\begin{definition}  A closed 1-form $\theta$ on $M$ is said to be
a {\em Morse 1-form} if $p^{*}(\theta) = df$ where $f$ is a Morse function on
$\widetilde{M}$.
\end{definition}

\begin{remarks} Another obviously equivalent way is to say that $\theta$ is a
Morse 1-form if locally (i.e. in a neighborhood of every given point in
$M$) the form
$\theta$ can be presented as $\theta = df$  where $f$ is a Morse function in
this neighborhood.
In fact it is sufficient to have such a presentation in a neighborhood of any
zero of $\theta$.
\end{remarks}

Assume for the rest of this section that $\theta$ is a Morse 1-form and that
$p^{*}(\theta)=df$ where $f$ is a Morse function on $\widetilde{M}$.

\begin{definition}  The {\em index} of a zero point $x$ of $\theta$ is
by fiat the index of a critical point $\stackrel{\sim}{x}$ of $f$ for any
lift $\stackrel{\sim}{x}$ of $x$.
\end{definition}

The following local description for Morse 1-forms $\theta$ is just a
reformulation of the classical Morse lemma (cf. \cite{Mi})
and can be deduced from the equivalent definition of a Morse 1-form given
in the remarks 3.2 above.

\begin{lemma}[Morse lemma for Morse 1-forms] Let $\theta$ be a Morse 1-form
on $M$.

\begin{itemize}
\item[(a)] Suppose that $p$ is a point such that
$\theta_{p}\neq 0$.  Then there is a chart $({\mathcal U},\phi)$ centered at
$p$ such that
$$
\left(\phi^{-1}\right)^{*}\theta\left(u_{1},\ldots,u_{n}\right) = du_{1}.
$$

\item[(b)] Suppose that $p$ is a zero point of $\theta$ of
index $k$.  Then there is a chart $({\mathcal U},\phi)$ centered at $p$
such that
$$
\left(\phi^{-1}\right)^{*}\theta\left(u_{1},\ldots,u_{n}\right) =
-\sum^{k}_
{j=1}u_{j}du_{j}+\sum_{j=k+1}^{n}u_{j}du_{j}.
$$
\end{itemize}
\end{lemma}

\begin{definition}.  Let $m_{k}(\theta)$ denote the number of zero
points of the Morse 1-form $\theta$ of index $k$. It is called the $k$-th
{\em Morse number} of $\theta$.
\end{definition}

The following is the analogue of the classical result which states that
on a compact closed manifold the set of Morse functions is open and dense
in the space of all smooth functions.

\begin{proposition}[Density of Morse 1-forms]
Let $M$ be a compact closed manifold.
 Then for every (de Rham) cohomology class $\alpha\in H^{1}(M;{\mathbb R} )$
the set of  Morse 1-forms $\theta\in\alpha$ is open and dense in the set of
all the forms in this cohomology class.
\end{proposition}

This proposition can be easily deduced from its classical analogue.

\subsection{Semiclassical asymptotics}
We will briefly describe an $L^2$-version of semiclassical asymptotics,
which are
similar to the ones which appear when we take the Witten deformation of the
de Rham
complex and consider the corresponding Laplacian (cf.
\cite{{Wi},{CFKS},{HS}}. For the
case of the algebra $\A$ corresponding to the regular representation of
$\pi_1(M)$,
such asymptotics were first proved in \cite{Sh}. The proofs given in
\cite{Sh} work
in a more general situation which we need now.

Let $H=-hA+B+h^{-1}V$ be a second order differential operator
 acting in
sections of a  Hilbertian $\A$-bundle $\F\to M$ over $M$.
Let $F$ denote the Hilbertian $\A$-module which is the fiber of $\F$.
We assume that all the coefficients in the local representations of $H$ are
smooth functions with values in $\End_{\A}(F)$, so $H$ commutes with the
fiberwise
action of $\A$.

We will also assume that $H$ satisfies  requirements similar to the ones in
\cite{Sh}.
Namely, we require that
 $A$ is a  second order elliptic operator with a negative principal symbol,
$B$ is a zeroth order operator, $V$ is a non-negative potential function
on $M$ which has only nondegenerate zeroes and $h>0$ is a small  parameter.

 Let   $N_\A(\lambda;\,H)$ denote the von Neumann
spectral density function
of the operator $H$. Let $K$ denote the {\em model operator} of $H$ (cf.
\cite{Sh}),
which is obtained as a direct sum of quadratic parts of $H$ in all zeros of $V$.
More precisely, let $\{\bar x_j,j=1,\dots,N\}$ be the set of all zeros of
$V$. Then
$K = \sum_j K_{j}$, where
$$
K_{j}: L^2({\mathbb R}^n,F) \rightarrow L^2({\mathbb R}^n,F)
$$
corresponds to the zero $\bar x_j$, acts in the Hilbert space $L^2({\mathbb
R}^n,F)$ of
all $F$-valued $L^2$-functions on $\R^n$ and has the form
$$
K_{j}=-A_{j}^{(2)}+B_{j}+V^{(2)}_{j},
$$
where all the components are obtained from $H$ as follows. Let us fix local
coordinates
on $M$ and a flat trivialization of $\F$ near $\bar x_j$.
The second order term $A_{j}^{(2)}$ is a homogeneous second order
differential operator
with constant coefficients (without lower order terms) obtained by
isolating the second order
terms in the operator $A$ and freezing the coefficients of this operator at
$\bar x_j$. The zeroth order term $B_{j}$ is a constant $\A$-endomorphism of $F$
which is obtained by  freezing the coefficients
of  $B$ at  $\bar x_j$. The other zeroth
order term $V^{(2)}_{j}$ is obtained by
taking the quadratic part of $V$ in the chosen coordinates near
$\bar x_j$.

Let $\Big\{\mu_j: j=1,2,3....\Big\}$ be the  eigenvalues
of the model
operator $K$, $\mu_i\not=\mu_j$ for $i\not=j$, and $r_j$ denote the
multiplicity of $\mu_j$. Let $\sigma(H)$ denote the $L^2$-spectrum of $H$,
i.e. its spectrum in the Hilbert space of $L^2$-sections of $\F$.
Then the following result is an easy generalization
of the corresponding result in \cite{Sh}:

\begin{theorem}[Semiclassical Approximation]  For any $R>0$ and $\ka\in (0,1)$
there exist $C>0$  and $h_0>0$
 such that
$$\sigma(H)\cap[-R,R]\subset\cup_{j=1}^\infty
(\mu_j-Ch^{1/5},\mu_j+Ch^\ka)\,.$$
Moreover for any
$j=1,2,3,\dots$ {\it with} $\mu_j\in [-R,R]$ and any $h\in (0,h_0)$
one has
$$N_\A(\mu_j+Ch^\ka;\,H) - N_\A(\mu_j-Ch^{1/5};\,H)
=r_j\;.$$
Besides $h^\ka$  can be replaced by $\exp(-C^{-1}h^{-1+\eps})$
 and $h^{1/5}$  by $h^\ka$   if
$H$  is flat  near all points $\bx_j$ (i.e. if $H$ coincides with $K_{j}$
near $\bar x_j$
for every $j$).
\end{theorem}

This means that for small $h$ the spectrum of $H$ concentrates near the
eigenvalues of
the model operator $K$, and for every such eigenvalue, the von Neumann
dimension
of the spectral subspace
of the operator $H$, corresponding to the part of the spectra near the
eigenvalue, is exactly equal to the multiplicity of
this eigenvalue.

We now state the main result of this section.

\begin{theorem}  Let $M$ be a compact manifold, $f$ a Morse function on $M$
and $\theta$ a closed 1-form on $M$. Let $\E\to M$ be a flat Hilbertian
$\A$-bundle over $M$ with fiber $E$ which is a finitely generated
Hilbertian $\A$-module.

\begin{itemize}
\item[(1)] (Strong Morse Inequalities).
  $$\displaystyle (\dimt
E)^{-1}\sum_{j=0}^{k}(-1)^{k-j}b^{j}_{(2)}(M,\E_{\theta})
\leq\sum_{j=0}^{k}(-1)^{k-j}m_{j}(f)$$
for $k=0,1,2\ldots,n$, with equality when $k=\dim_{\R} M=n$.

\item[(2)]  (Asymptotic Strong Morse Inequalities). Assume additionally that
$\theta$ is a Morse form. Then for $s\gg 0$, one has
$$\displaystyle(\dimt E)^{-1}\sum_{j=0}^{k}(-1)^{k-j} b^{j}_{(2)}
({M},\E_{s\theta})
\leq \sum_{j=0}^{k}(-1)^{k-j}m_{j}(\theta)$$
for $k=0,1,2\ldots,n$, with equality when $k=n$.
\end{itemize}
\end{theorem}

Let $\alpha,\beta$ be closed  1-forms on $M$, so that $\alpha$
 is a Morse form. Note that $\beta+s\alpha$ will be a closed Morse 1-form
if $s$ is sufficiently large.  Consider the associated complex
$$
\left(\Omega_{(2)}^\bullet({M},\E),\nabla_{\beta+s\alpha}\right)
$$
The corresponding Laplacian
$\Delta_{\beta+s\alpha,j}=\nabla^{*}_{\beta+s\alpha}\nabla_{\beta+s\alpha}
+\nabla_{\beta+s\alpha}\nabla^{*}_{\beta+s\alpha}$
acts on  $\Omega^{j}_{(2)}({M},\E)$.

Define the Riemannian dual $V$ of $\alpha$ which is a vector field on $M$
satisfying $g(V,\cdot)=\alpha(\cdot)$. Let
$${\mathcal L}_V:\Om^\bullet(M,\E)\to\Om^\bullet(M,\E)$$
denote the Lie
derivative acting in $\E$-valued forms on $M$, i.e.
$${\mathcal L}_V(f\otimes\om)=\nabla f\otimes\om+f\otimes{\mathcal L}_V\om,\ \
f\in C^\infty(M,\E),\ \om\in\Om^\bullet(M)\;,$$
where ${\mathcal L}_V\om$
is defined as the usual Lie derivative, applied to a scalar form $\om$,
$\nabla$ is the
flat connection on $\E$.

\begin{lemma} $\;\Delta_{\beta+s\alpha,j}
=\Delta_{\beta,j}+s\left({\mathcal L}_{V}
+{\mathcal L}^{*}_{V}+2\langle\alpha,\beta\rangle\right)
+s^{2}|\alpha|^{2}$ where
${\mathcal L}_{V}$ denotes the Lie derivative of the vector field $V$,
defined as above,
and ${\mathcal L}_{V}^{*}$ denotes its $L^{2}$ adjoint.
\end{lemma}

\begin{proof} Since $(e(\alpha))^\ast=i(V)$, we obtain
$\;\Delta_{\beta+s\alpha,j}=\left\{\nabla_{\beta}+se(\alpha),\nabla^{*}_
{\beta}+si(V)\right\}_{+}$ where $\left\{A, B\right\}_{+} = AB + BA$
denotes the
anticommutator. Therefore
$$
\begin{array}{l}
\displaystyle
\Delta_{\beta+s\alpha,j}=\Delta_{\beta,j}+s\left\{\nabla_{\beta},i(V)
\right\}_{+}+
s\left\{e(\alpha),\nabla^{*}_{\beta}\right\}_{+}+s^{2}\left\{e(\alpha),i(V)
\right\}_{+} \\[+10pt]
\displaystyle
=\Delta_{\beta,j}+s({\mathcal L}_{V}+{\mathcal L}^{*}_{V}+2\langle
\alpha,\beta\rangle)+s^{2}|\alpha|^{2}\;.
\end{array}
$$
\end{proof}

%Note that the statement of Theorem 3.8 does not involve the metrics on
%$M$ and on $\E$. Therefore we are free to make judicious
%choices of metrics, which will be chosen to be {\em flat} near the zero
%points of the Morse 1-form $\alpha$.
Clearly
$$\;\frac{1}{s}\;\Delta_{\beta+s\alpha,j}=\frac{1}{s}\;
\Delta_{\beta,j}+({\mathcal L}_{V}+{\mathcal L}^{*}_{V}+2\langle
\alpha,\beta\rangle) +s|\alpha|^{2}$$
is of the form $H=-hA+B+h^{-1}V$
required to apply Theorem 3.7, where $h=\frac{1}{s}$, $A=-\Delta_{\beta,j}$ is
independent of $h$, $B={\mathcal L}_{V}+{\mathcal L}^{*}_{V}+2\langle
\alpha,\beta\rangle$ is a zeroth order operator and
$V=|\alpha|^{2}$ is a non-negative function with non-degenerate zeros
precisely at the
zeros of $\alpha$. Also all
the terms commute with the
action of $\A$.
Denote by
$$
K_{(j)}: \Om^j_{(2)}({\mathbb R}^n,E) \rightarrow \Om^j_{(2)}({\mathbb R}^n,E)
$$
the corresponding model operator. (Here $\Om^j_{(2)}({\mathbb R}^n,E)$ is
the Hilbert
space of $E$-valued $L^2$-forms of degree $j$ on $\R^n$.)
For $s$ large the operator $\frac{1}{s}\;\Delta_{\beta+s\alpha,j}$ develops
gaps in its spectrum which persist for all large values of $s$.  In the
semiclassical limit, the spectrum of
 $\frac{1}{s}\;\Delta_{\beta+s\alpha,j}$
converges to the spectrum of a model operator $K_{(j)}$,
which is the direct sum of
harmonic oscillator type operators which operate near the zero points of
$\alpha$.

More precisely, if $E^{j}_{\lambda}(\beta+s\alpha)$ denotes the spectral
projection of the operator $\frac{1}{s}\;\Delta_{\beta+s\alpha,j}$, then
it follows from Theorem 3.7 that the von Neumann dimension of the spectral
subspace $\Im\Big(E^{j}_{\varepsilon}(\beta+s\alpha)\Big)$ is a constant
which equals
$\dimt(\Ker K_{(j)})$ for all $s>\Big(\frac{C}{\varepsilon}\Big)^5$.
Note in particular that $\varepsilon$  lies in the gap of the spectrum
 of $\frac{1}{s}\;\Delta_{\beta+s\alpha,j}$ when $s\gg 0$.
The problem of calculating the eigenvalues (and their multiplicities)
of the model operator
$K_{(j)}$ actually reduces to the diagonalizing of $B$ at the critical
points of the Morse 1-form $\alpha$. Note that this is a local calculation,
which
was performed by Witten (cf. \cite{{Wi},{HS}}) in the standard case. An
obvious modification
of Witten's arguments to the case of forms with coefficients in $\E$
leads to the following

\begin{lemma}
The multiplicity of $0$  as an eigenvalue of $K_{(j)}$ is
$$
\dimt(\Ker K_{(j)})=m_{j}(\alpha)\ \dimt E.
$$
\end{lemma}

Let us take $L^{j}=\Im E^{j}_{\varepsilon}(\beta+s\alpha)$
for $s\gg 0$. Then it can be proved as in \cite{Sh2} that
$L^j$ is a finitely generated Hilbertian $\A$-module, and
$$
\nabla_{\beta+s\alpha}:L^{j}\rightarrow L^{j+1}
$$
is a bounded coboundary operator commuting with the $\A$-action;
in particular, $\nabla_{\beta+s\alpha}^2 = 0$.
The Hilbertian $\A$-complex
$$
L^\bullet :0\longrightarrow L^{0}\stackrel{\nabla_{\beta+s\alpha}}
{\longrightarrow} L^{1}
\stackrel{\nabla_{\beta+s\alpha}}{\longrightarrow}\ldots
\stackrel{\nabla_{\beta+s\alpha}}
{\longrightarrow}L^{n}\longrightarrow 0
$$
is  bounded chain homotopy equivalent to
$\left(\Omega_{(2)}^\bullet({M},\E),\nabla_{\beta+s\alpha}\right)$,
and in particular, its reduced $L^{2}$-cohomology
is $H^{j}_{(2)}(L^{\bullet},\nabla_{\beta+s\alpha})=
H^{j}_{(2)}({M}, \E_{\beta+s\alpha})$.

Now modifying an argument given in \cite{Sh2} we obtain

\begin{proposition} Let $M$ be a compact closed manifold,
 $\alpha,\beta$  closed  1-forms on $M$, and $\alpha$ is a Morse form.
Let $\E\to M$ be a flat Hilbertian $\A$-bundle
over $M$. Then the twisted de Rham complex of $L^2$ differential forms
$\left(\Omega_{(2)}^\bullet({M},\E),\nabla_{\beta+s\alpha}\right) $
is bounded chain homotopy equivalent to the finitely generated
Hilbertian $\A$-complex $\left(L^\bullet, \nabla_{\beta+s\alpha}\right)$.
\end{proposition}

\begin{proof}[Proof of Theorem 3.8 ]
By Lemmas 1.1 and 1.2, it follows that for $s\gg 0$ one
has
$$
(\dimt E)^{-1}\sum_{j=0}^{k}(-1)^{k-j}b^{j}_{(2)}({M},
\E_{\beta+s\alpha})\leq\sum_
{j=0}^{k}(-1)^{k-j}m_{j}(\alpha)
$$
with the equality for $k=n$.

Setting $\alpha=df$ and $\beta=\theta$,
we observe that $[\beta+s\alpha]=[\theta]$ for all $s$,
proving  part (1).

Setting $\alpha=\theta$ and $\beta=0$, we observe that
$[\beta+s\alpha]=s
[\alpha]$, proving  part (2).
\end{proof}

The equalities (the case $k=n$) in part (2) of Theorem 3.8 follow also from
the stability of the
($L^2$)-index under deformations, so in fact it is not necessary to take
$s\gg 0$.
More precisely, one has

\begin{proposition} For any closed Morse 1-form $\theta$ on $M$ one has
$$\displaystyle (\dimt E)^{-1}\sum_{j=0}^{n}(-1)^{j}b^{j}_{(2)}({M},\E_\theta)
 ={\chi}(M)= \sum_{j=0}^{n}(-1)^{j}m_j(\theta).$$
\end{proposition}

\begin{proof}  The first equality is  in fact a particular case of the general
$L^2$-index theorem by I.Singer \cite{Si}. (See also \cite{A} for the case
of the algebra $\A$ generated by the regular representation of $\pi_1(M)$.)
Firstly, observe that for the
$L^{2}$-index of the operator
$\nabla_{\theta} +\nabla^{*}_{\theta}$ we have
$$
\begin{array}{lcl}
\displaystyle\mbox{index}_{L^{2}}(\nabla_{\theta}+\nabla^{*}_{\theta}) &=&
\dimt
\Ker(\nabla_{\theta}+\nabla^{*}_{\theta})|_{\Omega_{(2)}^{even}({M},\E)} -
\dimt \Ker(\nabla_{\theta}+\nabla^{*}_{\theta})|_{\Omega_{(2)}^{odd}({M},\E)}\\
& = &\sum_{j=0}^{n}(-1)^{j}b^{j}_{(2)}({M},\E_\theta),
\end{array}
$$
where $\Omega_{(2)}^{even}(M,\E)$ denotes the even degree $L^2$ differential
forms on $M$ with values in $\E$, and $\Omega_{(2)}^{odd}(M,\E)$ the odd
degree ones. By the deformation invariance of the $L^{2}$-index, one sees that
$\mbox{index}_{L^{2}}(\nabla_{\theta}+\nabla^{*}_{\theta})=
\mbox{index}_{L^{2}}(\nabla+\nabla^{\ast})$, where $\nabla+\nabla^{\ast}$
is considered as an operator
in $\Omega_{(2)}^{\bullet}({M},\E)$. Together with Theorem 3.8,
one deduces the proposition.
\end{proof}

The idea of the following proposition is due to M. Gromov.
 A quantitative version of this proposition will be given in the next section.

\begin{proposition} Let $\theta$ be a closed 1-form on $M$, and
$\sigma(\Delta_{\theta,j})$ denote the $L^2$-spectrum of the twisted
Laplacian $\Delta_{\theta,j}$.

\begin{itemize}
\item[(a)] Suppose that ${0}\in\sigma(\Delta_{\theta,j})$.
Then $m_{j}(f)>0$ for every Morse function $f$ on $M$.

\item[(b)] Suppose that $\theta$ is a Morse form and
${0}\in\sigma(\Delta_{s\theta,j})$ for $s\gg 0$.  Then
$m_{j}(\theta)>0$.
\end{itemize}
\end{proposition}

\begin{proof}
(a) The key idea is to use the abstract Theorem 1.3 (see also \cite{GS})
which shows that $0$-in-the-spectrum is a
homotopy-invariant phenomenon.  From this theorem  we deduce first
that the inclusions ${0}\in\sigma(\Delta_{\theta+sdf,j})$,
depending formally on $s$, are equivalent for different $s$ because the
corresponding
complexes are isomorphic (see  proof of Proposition 2.1). Secondly, for
$s\gg 0$ the twisted complex $(\Om^\bullet_{(2)}(M,\E),\nabla_{\theta+sdf})$ is
homotopy equivalent to the complex $(L^\bullet,\nabla_{\theta+sdf})$ with
$\dimt L^j=m_j(f)\dimt E$  due to Lemma 3.10 and Proposition 3.11.
Since $m_j(f)=0$ would imply that
$0$ is not in the spectrum of the Laplacian in the complex
$(L^\bullet,\nabla_{\theta+sdf})$, the same is true in
 the homotopy equivalent
complex $(\Om^\bullet_{(2)}(M,\E),\nabla_{\theta+sdf})$, and (a) follows.

The proof of (b) can be done by similar arguments.
\end{proof}

\section{Asymptotic $L^2$ Morse-Farber inequalities}
In this section
we briefly review the extended category construction
of Farber \cite{F} and we use it to define the
the extended twisted de Rham $L^2$-cohomology. (See also
L\"uck \cite{Lu2} for a different and purely algebraic approach
to $L^2$-cohomology which is  equivalent to the Farber's approach.)
The main result of the section is the
asymptotic $L^2$ Morse-Farber inequalities
for the extended twisted de Rham $L^2$-cohomology.
We shall assume in this section that $\A$ is
a finite von Neumann algebra equipped
with a finite, normal and faithful trace $\tau$.

\subsection{The extended category}
In \cite{F} Farber used a P.Freyd construction to define  an Abelian category
$\E(\A)$, which he called
the extended category of Hilbertian $\A$ modules.
It contains the usual
additive category $\H(\A)$ of Hilbertian
$\A$ modules,  which turns out to
 be the full subcategory of projective objects in
$\E(\A)$. We refer to \cite{F} for further details
on the following discussion about $\E(\A)$.

 An object  in
$\E(\A)$ is a morphism $(\alpha:A'\to A)$ in $\H(\A)$
(so $A,A'$ are objects in $\H(\A)$), and it is
called a {\em virtual Hilbertian $\A$ module}.

A morphism
$$(\alpha:A'\to A)\longrightarrow(\beta:B'\to B)$$
in $\E(\A)$ is defined by a morphism $f:A\to B$ in $\H(\A)$ such that there
exists
a morphism $h:A'\to B'$ (in $\H(\A)$) with $f\alpha=\beta h$. Two such
morphisms
$f,f':A\to B$  define the same morphism in $\E(\A)$ if and only if
there exists a morphism $g:A\to B'$ such that $f-f'=\beta g$.

A Hilbertian $\A$ module $A$ can be canonically viewed
as a virtual Hilbertian $\A$ module $(0\to A)$, and
this correspondence is an embedding of $\H(\A)$ into $\E(\A)$.

A virtually Hilbertian $\A$-module $(\alpha: A'\to A)$
is said to be a {\em torsion} object in $\E(\A)$ if
$\overline{\Im \alpha} = A$.
The Novikov-Shubin invariant is an
invariant of a torsion object in $\E(\A)$.
A virtually Hilbertian $\A$ module $(\alpha: A'\to A)$
is  a {\em projective} object in $\E(\A)$ if and only if
 $\Im\alpha$ is closed.

The von Neumann dimension is an invariant of a projective
object in $\E(\A)$.

A virtually Hilbertian $\A$ module $(\alpha: A'\to A)$ has
a canonically defined torsion part $(\alpha: A'\to \overline{\Im \alpha})$
and a canonically defined projective part
$(0\to A/\overline{\Im \alpha})$, and it is a direct sum of these two parts
(though not canonically).

\subsection{The extended $L^2$-cohomology}
Let
$$C^\bullet:\quad \dots \to C^{i-1}\to C^i\to C^{i+1}\to \dots$$
be a Hilbertian $\A$-complex (cf. section 1). It can be
viewed as a complex of virtually Hilbertian $\A$-modules,
i.e. as a complex in the extended category $\E(\A)$.
Then the $i$-th  cohomology of  $C^\bullet$ is well defined in $\E(\A)$
(because $\E(\A)$ is an Abelian category), and is explicitly
given as follows:
$$
\h^i(C^\bullet) = (d : C^{i-1} \to Z^i),
$$
where $Z^i$ denotes the Hilbertian $\A$-submodule
of cocycles in $C^i$. It is called the $i$-th
{\em extended
$L^2$-cohomology} of the Hilbertian $\A$-complex
$C^\bullet$. The projective part of the extended
$L^2$-cohomology is the Hilbertian $\A$-module
$$
P(\h^i(C^\bullet)) = Z^i/ \overline{\Im( d: C^{i-1}\to C^i)} =
H^i_{(2)}(C^\bullet)
$$
which coincides with the reduced $L^2$-cohomology of
the Hilbertian $\A$-complex $C^\bullet$.
The torsion part of the extended
$L^2$-cohomology is the morphism
$$
T(\h^i(C^\bullet)) = (d: C^{i-1}\to \overline{\Im( d: C^{i-1}\to C^i)})
$$

\subsection{The extended twisted de Rham $L^2$-cohomology}
Let $M$ be a compact manifold
and $\beta$ a closed 1-form on $M$. Let $\E\to M$ be a flat Hilbertian
$\A$-bundle
over $M$. Following \cite{Sh2} let us define the
{\em extended twisted de Rham $L^2$-cohomology}
of the complex
$$
\left(\Omega_{(2)}^\bullet({M},\E),\nabla_{\beta}\right)
$$
as the extended cohomology of any finitely
generated Hilbertian $\A$-complex which is bounded chain
homotopy equivalent to the given de Rham complex
$\left(\Omega_{(2)}^\bullet({M},\E),\nabla_{\beta}\right)$.
Since bounded chain homotopy
is an equivalence relation, this is well defined.
For example, let $E^{j}_{\lambda}(\beta)$ be
the spectral
projection of the operator $\Delta_{\beta,j}$. Let us fix $\eps>0$ and define
$L^{j}=E^{j}_{\varepsilon}(\beta)$. Then $L^j$
can be shown to be a finitely generated Hilbertian $\A$-module
(cf. \cite{Sh2}), and
$$
\nabla_{\beta}:L^{j}\rightarrow L^{j+1}
$$
is the bounded coboundary operator commuting with the $\A$-action.
In particular, $\nabla_{\beta}^2 = 0$.
The Hilbertian $\A$-complex
$\left(L^\bullet, \nabla_{\beta}\right)$
was earlier observed to be bounded chain homotopy
equivalent to
$\left(\Omega_{(2)}^\bullet({M},\E),\nabla_{\beta}\right)$.
By definition, the {\em extended twisted  de Rham $L^2$-cohomology}
is the extended cohomology of the finitely
generated Hilbertian $\A$-complex
$\left(L^\bullet, \nabla_{\beta}\right).$
It is denoted by $\h^\bullet(M,\E_\beta)$ and is
represented in $\E(\A)$ by the morphism
$$
\h^i(M,\E_\beta) = (\nabla_{\beta} : L^{i-1} \to Z^i)
$$
where $Z^i$ denotes the Hilbertian $\A$-submodule
of cocycles in $L^i$.
The projective part of the extended twisted
de Rham $L^2$-cohomology is the Hilbertian $\A$-module
$$
P(\h^i(M,\E_\beta)) = Z^i/ \overline{\Im(\nabla_{\beta}: L^{i-1}\to L^i)}
= H^{i}_{(2)}(L^{\bullet},\nabla_{\beta})
$$
which coincides with the reduced
de Rham $L^2$-cohomology of
the Hilbertian $\A$-complex $\left(L^\bullet, \nabla_{\beta}\right).$
But as observed earlier, the reduced $L^{2}$-cohomology
$H^{i}_{(2)}(L^{\bullet},\nabla_{\beta})$ coincides with
the reduced twisted  de Rham $L^2$-cohomology
$H^{i}_{(2)}({M}, \E_{\beta})$.
The torsion part of the extended twisted
de Rham $L^2$-cohomology  is the morphism
$$
T(\h^i(M,\E_\beta)) = (\nabla_{\beta}: L^{i-1}\to
\overline{\Im( \nabla_{\beta}: L^{i-1}\to L^i)})\;.
$$

Alternatively we can use Witten deformation and take
$L^j_W=\tilde E^j_\eps(\beta+sdf)$ instead of $L^j$. Here
$f$ is a Morse function on $M$, $\tilde E^j_\eps(\beta+sdf)$
is the spectral projection of $\frac{1}{s}\De_{\beta+sdf,j}$,
$\eps>0$ is fixed, $s\gg 0$. Then  again $L^j_W$ will be a finitely generated
Hilbertian $\A$-module (\cite{Sh2}), and
$$\nabla_{\beta+sdf}:L^j_W\to L^{j+1}_W$$
is a bounded coboundary operator, defining a Hilbertian $\A$-complex
$(L^\bullet_W,\nabla_{\beta+sdf})$ which is
bounded chain homotopy equivalent to $(\Om^\bullet_{(2)},\nabla_\beta)$.
The advantage of this approximation is that $\eps$ will be in the gap of
the spectrum
of the Laplacian $\frac{1}{s}\De_{\beta+sdf}$ for $s\gg 0$, so the spectral
subspace
$L^j_W$ has better analytic and geometric properties.

\subsection{Minimal number of generators}
M. Farber \cite{F} introduced the minimal number of generators as an invariant
of virtual Hilbertian $\A$-modules. It is  non-trivial
on all non-trivial torsion modules and takes
 values in the nonnegative integers. More precisely,
let $\chi$ be a virtually Hilbertian $\A$-module, and
$\mu(\chi)$ be the smallest integer $\mu$ such that there exists
an epimorphism from the direct sum of $\mu$ copies of $\ell^2(\A)$
onto $\chi$. Then $\mu(\chi)$ is called the {\em minimal number of
generators} of $\chi$ and it has the following properties:

\begin{itemize}
\item[(1)] $\mu(\chi) = 0$ if and only if $\chi=0$; besides,
$\mu(\chi) = \mu(P(\chi))$ if and only if $T(\chi) = 0$.

\item[(2)] If $\chi$ is projective, then $\mu(\chi) \ge \dimt \chi$.

\item[(3)] if $\chi'$ is another virtually Hilbertian $\A$-module, then
$$
\max\{\mu(\chi), \mu(\chi')\} \le \mu(\chi\oplus\chi') \le \mu(\chi) +
\mu(\chi').
$$

\end{itemize}

Some calculations of this invariant are done in \cite{F}, section 7.

The minimal number of generators seems to be useful for torsion objects
provided $\A$
has a non-trivial center. It is easy to see that   if $\A$ is a factor,
then $\mu(\chi)=1$
for any non-trivial torsion object.

Indeed, if $\A$ is a factor, then any
torsion module $\chi$  can be represented by  a morphism
$\alpha: M^\prime\to M$
where $\dim_\A M$ is arbitrarily small. (This follows from Farber's
excision property \cite{F}
if we apply a  spectral cut to $\alpha$, removing a part of the spectrum
of $|\alpha|$ outside $[0,\epsilon]$ for small $\epsilon$.) It follows that
there is an epimorphism from $M$ (identified with $0\to M$) onto $\chi$.
Now by the fundamental property  of Hilbert modules over factors
(see e.g. \cite{Di} or \cite{T})
we can embed $M$ into $l^2(\A)$ as a Hilbert $\A$-module,
and, therefore,  produce epimorphism (by orthogonal projection) of
$l^2(\A)$ onto $M$,
hence onto $\chi$. This means that  $\mu(\chi)\le 1$.

The following abstract theorem is due to Farber \cite{F},
theorem 8.1. It has been rephrased here in terms of cohomology.

\begin{theorem} Let
$$C^\bullet :\quad \dots \to C^{i-1}\stackrel{d}{\rightarrow} C^i
\stackrel{d}{\rightarrow}C^{i+1}\to \dots$$
be a free finitely generated Hilbertian chain complex in $\E(\A)$.
Then for any integer $i$ the following inequality holds:
$$
\dimt (C^i)\ \ge \ \mu[\h^i(C)\oplus T(\h^{i+1}(C))].
$$
\end{theorem}

The first part of the following Theorem was essentially proved by Farber
\cite{F}
using combinatorial methods. Our method is different (we use analysis),
and also the second part of the Theorem is new.
\begin{theorem}  Let $M$ be a compact manifold, $f$ a Morse function on $M$
and $\theta$ a Morse 1-form on $M$. Let $\E\to M$ be a flat Hilbertian
$\A$-bundle
over $M$ such that its fiber $E$ is a  free finitely generated
Hilbertian $\A$-module.
%Assume also that $\A$ is a factor.

\begin{itemize}
\item[(i)] ($L^2$ Morse-Farber Inequalities).
$$
m_j(f)\ \ge\ (\dimt E)^{-1}\cdot
\mu[\h^j(M,\E_\theta)\oplus T(\h^{j+1}(M,\E_\theta))]
$$
for $j=0,1,2,\dots$.

\item[(ii)]  (Asymptotic $L^2$ Morse-Farber Inequalities).
For $s\gg 0$, one has
$$
m_j(\theta)\ \ge\ (\dimt E)^{-1}\cdot
\mu[\h^j(M,\E_{s\theta})\oplus T(\h^{j+1}(M,\E_{s\theta}))]
$$
for $j=0,1,2,\dots$.
\end{itemize}
\end{theorem}

\begin{proof} Let $\alpha,\beta$ be closed 1-forms on $M$, and besides $\alpha$
is a Morse form. It is convenient to use the Witten ``small eigenvalues"
approximation
$(L^{\bullet}_W,\nabla_{\beta+s\alpha}),\ s\gg 0,$ of the complex
$(\Om^\bullet_{(2)}(M,\E_\beta),\nabla_\beta)$
as described at the end of sect. 4.3.

Assume first that  $\A$ is a factor. Since $E$ is a free
Hilbertian $\A$ module, it follows that $L^j_W$ is a
free Hilbertian $\A$-module for all $j$, because $\dimt L^j_W$ is an integer
by Lemma 3.10 and the discussion before it. For general $\A$ we still can prove
that $L^j_W$ will be free by taking the factor decomposition of $E$ and
using arguments from \cite{Sh2}.

An alternative way to prove that $L^j_W$ is free: observe that a small
refinement of the arguments in \cite{Sh} shows that $L^j_W$ is isomorphic
(as a Hilbertian $\A$-module)
to $\Ker K_{(j)}$, where $K_{(j)}$ is the model operator corresponding to the
Hamiltonian $H=\frac{1}{s}\De_{\beta+s\alpha,j}$ as described in section 3.2.

The theorem follows if we  apply
Theorem 4.1 to the the free finitely generated
Hilbertian $\A$-complex
$\left(L^\bullet_W, \nabla_{\beta+s\alpha}\right)$, for $s\gg 0$,
as in Proposition 3.11.
Setting $\alpha=df$ and $\beta=\theta$,
we observe that
$[\beta+s\alpha]=[\theta]$ for all $s$, proving part (i).

Setting $\alpha=\theta$ and $\beta=0$, we observe that
$[\beta+s\alpha]=s[\alpha]$,
proving part (ii).
\end{proof}

The following corollary can be viewed as a quantitative version of
Proposition 3.13.

\begin{corollary} Let $M$ be a compact manifold, $f$ a Morse function on $M$
and $\theta$ a closed 1-form on $M$. Let $\E\to M$ be a flat Hilbertian
$\A$-bundle
over $M$ such that its  fiber is a free finitely generated
Hilbertian $\A$-module $E$.
Let $\lambda_0(\Delta_{\theta,j})$ denote the
bottom of the $L^2$-spectrum of the twisted Laplacian $\Delta_{\theta,j}$
acting on the complement of its $L^2$-kernel.
%Assume also that $\A$ is a factor.

\begin{itemize}
\item[(a)] Suppose that $\lambda_0(\Delta_{\theta,j}) = 0$, i.e. there is
no spectral
gap at zero.
Then
$$
m_j(f)\ >\ (\dimt E)^{-1}\cdot
\mu[P(\h^j(M,\E_\theta))] \ge (\dimt E)^{-1}\cdot b_{(2)}^j(M,\E_\theta).
$$
for every Morse function $f$ on $M$.

\item[(b)] Suppose  that $\theta$ is a Morse 1-form on $M$ and
$\lambda_0(\Delta_{s\theta,j}) = 0$
for $s\gg 0$.
 Then
$$
m_{j}(\theta) \ >\ (\dimt E)^{-1}\cdot
\mu[P(\h^j(M,\E_{s\theta}))] \ge (\dimt E)^{-1}\cdot b_{(2)}^j(M,\E_{s\theta}).
$$
for $s\gg 0$.
\end{itemize}
\end{corollary}

\begin{proof}
Since $\lambda_0(\Delta_{\theta,j}) = 0$, it follows that
$T(\h^{j}(M,\E_\theta))$ is a non-trivial virtual Hilbertian
$\A$-module. By property (1) of the minimal number of generators, one has
$\mu(\h^{j}(M,\E_\theta)) > \mu(P(\h^{j}(M,\E_\theta)))$.
By Theorem 4.2, part (i), and property (3) of the minimal number of generators,
one has
$$
m_j(f)\ \ge\ (\dimt E)^{-1}\cdot
\mu[\h^j(M,\E_\theta)] > (\dimt E)^{-1}\cdot
\mu[P(\h^j(M,\E_\theta))].
$$
The last inequality in part (a) follows from property (2) of the minimal
number of
generators.

Part (b) is proved similarly.
\end{proof}

\section{Calculations.}

In this section, we do some calculations in the special case of the Hilbertian
$({\mathcal U}(\pi)-\pi)$-bimodule $\ell^2(\pi)$. Let $\E\to M$ denote the
associated
flat Hilbertian ${\mathcal U}(\pi)$-bundle over $M$. Then it is well known that
the Hilbertian ${\mathcal U}(\pi)$-complexes
$\left(\Omega_{(2)}^\bullet({M},\E),\nabla\right)$ and
$\left(\Omega_{(2)}^\bullet(\widetilde{M}), d\right)$ are canonically
isomorphic, where $p:{\widetilde M}\to M$ denotes the universal
cover of $M$. This isomorphism also establishes that if $\theta$ is a closed
1-form on $M$, then
$\left(\Omega_{(2)}^\bullet({M},\E),\nabla_{\theta}\right)$ and
$\left(\Omega_{(2)}^\bullet(\widetilde{M}), d_{\theta}\right)$
are canonically isomorphic, where
$d_\theta = d + e(p^*\theta)$. Here $ e(p^*\theta)$ denotes exterior
multiplication by the closed 1-form $p^*\theta$. In this case, we shall
denote the twisted $L^2$-invariants by $R(\widetilde M, \theta) =
R(M,\E_\theta)$, where $R$ denotes any particular $L^2$-invariant as in the
previous sections.

We begin with the following basic lemma.

\begin{lemma}  Let $M$ be a closed connected manifold and $\pi_1(M)$ be
infinite.  Then
$$
H^{0}_{(2)}(\widetilde{M},\theta)=\{0\}=H^{n}_{(2)}(\widetilde{M},\theta),
$$
where $n=\dim_{\R}M$.
\end{lemma}

\begin{proof}  Since
$H^{0}_{(2)}(\widetilde{M},\theta)=\Ker \De_{\theta,0}=\Ker d_{\theta}$ on
$\Om_{(2)}^0(\widetilde M)$,
we will first prove that $\Ker d_{\theta}=\{0\}$ on
$L^{2}$-functions
on $\widetilde{M}$.  Let $g\in\Ker d_{\theta}$.  Let $p^\ast\theta=df$,
$f\in C^\infty(\widetilde M)$. Since
$$
0=d_{\theta}g=e^{-f}d(e^{f}g),
$$
we see that $d(e^{f}g)=0$, hence $e^{f}g=c=\operatorname{const}$
and also $g\in L^{2}(\widetilde{M})$.

Now let us check that the inclusion $g=ce^{-f}\in L^2({\widetilde M})$ is
possible
only if $c=0$. Let $F$ be a fundamental domain of the action of
$\pi=\pi_1(M)$ on
${\widetilde M}$ by deck transformations, so that
$${\widetilde M}=\cup\{\ga F|\;\ga\in\pi\}$$
up to a set of measure $0$. We can assume $F$ to be open and connected.
Let us fix a point $P\in F$ and assume that $f(P)=0$. Then $f(\ga
P)=\int_{l_\ga}\theta$ where
$l_\ga$ is a loop representing $\ga$ in $\pi_1(M)$ with the base point $P$.
Clearly
$f(\ga_1\ga_2P)=f(\ga_1P)+f(\ga_2P)$. Denote $\phi(\ga)=\exp(f(\ga P))$. Then
$\phi(\ga_1\ga_2)=\phi(\ga_1)\phi(\ga_2)$, so $\phi:\pi\to (0,\infty)$ is a
group homomorphism.

It is easy to see that for any $Q\in F$
$$f(\ga Q)-f(Q)=\int_{l_\ga}\theta=f(\ga P)-f(P)=f(\ga P)\,,$$
hence
$$\exp(-f(\ga Q))=\phi(\ga^{-1})\exp(-f(Q)), \ \ Q\in F,\ \ga\in\pi\;.$$
Integrating over $F$ and summing over all $\ga\in\pi$, we see that
$$\int_{\tilde M}|\exp(-f(x))|^2dx=
\left(\sum_{\ga\in\pi}|\phi(\ga)|^2\right)\int_F |\exp(-f(x))|^2dx\,,$$
where $dx$ denotes any $\pi$-invariant smooth measure with a positive
density on
${\widetilde M}$. It follows that the inclusion $e^{-f}\in L^2({\widetilde
M})$ is equivalent
to $\sum_{\ga\in\pi}|\phi(\ga)|^2<\infty$ which is obviously impossible for
infinite $\pi$.
Therefore the constant $c$ above should vanish and $g\equiv 0$.

We conclude that $H^{0}_{(2)}(\widetilde{M},\theta)=\{0\}$.

By Poincar\'{e} duality argument (applied on ${\widetilde M}$ which is
orientable!)
we obtain
$$H^{n}_{(2)}(\widetilde{M},\theta)=H^{0}_{(2)}(\widetilde{M},-\theta) =
\{0\}\;.$$
\end{proof}

Our next result proves the upper semi-continuity property of the twisted
$L^2$ Betti numbers as a function on $H^{1}(M,{\mathbb R} )$.

\begin{lemma}
The function
$H^{1}(M,{\mathbb R} )\to {\mathbb R} $
given by
$[\theta] \mapsto b^{j}_{(2)}(\widetilde{M}, \theta)$
is upper semi-continuous.
\end{lemma}

\begin{proof}
We identify $H^{1}(M,{\mathbb R} )\cong{\mathcal H}^{1}(M,{\mathbb R} )$.
The proof is based on the fact that
$\theta\mapsto\Delta_{\theta,j}$
is an analytic family.
Recall that a function $f:{\mathcal H}^{1}
(M,{\mathbb R} )\rightarrow {\mathbb R} $ is said to be upper semi-continuous at
$\alpha\in{\mathcal H}^{1}(M,{\mathbb R} )$ if for any $\varepsilon>0$,
there exists a
$\delta>0$ such that for any $\theta$ satisfying $|\theta-\alpha|<\delta$,
we have
$$
f(\theta)<f(\alpha)+\varepsilon
$$
For each $\alpha\in{\mathcal H}^{1}(M,{\mathbb R} )$,
the operator $\Delta_{\alpha,j}$
is $\pi$-Fredholm in $\Omega^{j}_{(2)}(\widetilde{M})$. (Here $\pi=\pi_1(M)$.)
This means firstly that  $\dim_\pi\Ker \Delta_{\alpha,j}<\infty$
where $\dim_\pi=\dimt$ for
the canonical trace $\tau$ on $\U(\pi)$, and secondly that
for any $\eps>0$ there exists a closed $\pi$-invariant subspace
$L_\eps\subset\Im\De_{\alpha,j}$ such that
$\dim_\pi(\overline{\Im\De_{\alpha,j}}\ominus L_\eps)<\eps$ i.e.
$\overline{\Im\De_{\alpha,j}}=L_\eps\oplus L^\prime_\eps$ with
$\dim_\pi L^\prime_\eps<\eps$.

Note also that
$$
\Omega^{j}_{(2)}(\widetilde{M})=\Ker\Delta_{\alpha,j}
\oplus\overline{\Im\Delta_{\alpha,j}}\,.
$$

By the closed graph theorem, there exists $C>0$ such that
$$
\|\Delta_{\alpha,j}u\|\geq C\|u\|,\ \ u\in L_{\varepsilon}\,.
$$
It follows that for any $\de^\prime>0$ we can choose $\de>0$ such that
$$
\|\Delta_{\theta,j}u\|\geq (C-\delta')\|u\|\;\;
\mbox{ whenever } u\in L_\eps\;\;\mbox{ and }\;\;|\theta-\alpha|<\delta.
$$
Therefore
$$
\dim_{\pi}(\Ker\Delta_{\theta,j})\leq
\mbox{codim}_{\pi}(L_{\varepsilon})
<\dim_\pi(\Ker\Delta_{\alpha,j})+\varepsilon.
$$
\end{proof}

Let ${\mathcal L}_{{\mathbb Q}}$ denote the class of finitely presented discrete
groups $\pi$ such that the von Neumann dimension of the kernel of any
$\pi$-invariant operator acting on $\ell^{2}(\pi)\otimes {\mathbb C}^l$
and which comes from the group algebra ${\mathbb C}(\pi)$, is a rational
number.
${\mathcal L}_{{\mathbb Q}}$ contains all elementary amenable
groups and extensions
of groups in ${\mathcal L}_{{\mathbb Q}}$ by right orderable groups.    Let
${\mathcal L}_
{{\mathbb Z}}\subset{\mathcal L}_{{\mathbb Q}}$ denote the class of
finitely presented
discrete groups $\pi$ such that the von Neumann dimension
of the kernel of any $\pi$-invariant operator acting on $\ell^2(\pi)\otimes
{\mathbb C}^l$
and which comes from the group algebra ${\mathbb C}(\pi)$, is an integer.
  It is known that ${\mathcal L}_{{\mathbb Z}}$ contains all
torsion-free
elementary amenable groups and torsion-free extensions of groups in
${\mathcal L}_{{\mathbb Z}}$ by torsion-free right orderable groups.
It has been conjectured that ${\mathcal L}_{{\mathbb Q}}$ contains all finitely
presented groups and ${\mathcal L}_{{\mathbb Z}}$ contains all torsion-free
finitely presented groups.
(See Cohen \cite{C}, Donnelly \cite{Don} and especially Linnell
\cite{Li} for further information on this discussion.)

\begin{corollary}
Suppose that $b^{j}_{(2)}(\widetilde{M}) =0$ and
$\pi_1(M)\in {\mathcal L}_{{\mathbb Z}}$. Then there is a $\delta>0$ such that
$b^{j}_{(2)}(\widetilde{M},[\theta]) = 0$ for all $[\theta]$ with
$\|[\theta]\|<\delta$.
\end{corollary}

\subsection{ Flat manifolds}
Let $\theta$ denote the $S^1$-invariant nowhere zero 1-form on the circle
$S^{1}=\R/\Z$, which represents a generator for $H^{1}(S^{1},{\mathbb R}
)$.  Then
any other element of $H^{1}(S^{1},{\mathbb R} )$ is represented by $s\theta$,
$s\in {\mathbb R} $. We can assume that $|\theta(x)|=1$ in the canonical metric
for all $x\in S^1$. Choosing the sign, we can take $\theta=dx$ where $x$ is
the canonical
coordinate in $\R$.

On ${\mathbb R} $, and on functions, one calculates for the case of the
standard flat metric
$$
\Delta_{s\theta,0}=-\frac{\partial^{2}}{\partial x^{2}}+s({\mathcal L}_{V}+
{\mathcal L}_{V}^{*})+s^{2}\,,
$$
where $V=\frac{\partial}{\partial x}$ is dual to $dx$.
So ${\mathcal L}_{V}=\frac{\partial}{\partial x}=-{\mathcal L}^{*}_{V}$.
That is,
$$
\Delta_{s\theta,0}=-\frac{\partial^{2}}{\partial x^{2}}+s^{2}\;.
$$
We get a spectral resolution of $\Delta_{s\theta,0}$ using the Fourier
transform, and one calculates that
$$
\sigma(\Delta_{s\theta,0})=[s^{2},\infty)=\sigma(\Delta_{s\theta,1}).
$$
We conclude that
$$
\begin{array}{lclcl}
\lambda_{0,0}(\Delta_{s\theta}) & = & \lambda_{0,1}
(\Delta_{s\theta}) & = &
s^{2} \\
 &&&& \\
\alpha_{0}(\widetilde{S}^1,s[\theta]) & = &
\alpha_{1}(\widetilde{S}^1,s[\theta]) &
= & \left\{\begin{array}{lcl}
\frac{1}{2} & \mbox{ if } & s=0 \\[+5pt]
\infty & \mbox{ if } & s\neq 0
\end{array}\right.      \\
 &&&& \\
b^{0}_{(2)}(\widetilde{S}^1,s[\theta]) & = &
b^{1}_{(2)}(\widetilde{S}^1,s[\theta])
& = & 0\;\mbox{ for all }\;s.
\end{array}
$$
Here $\lambda_{0,j}(A)$ denotes the bottom of the spectrum
of the operator $A$
acting on $L^2$-forms of degree $j$.
Generalizing this calculation to the higher dimensions, one has

\begin{proposition}.  Let $\theta$ be a closed 1-form on a compact flat
manifold $M$ of dimension $n$ (i.e. $M$ is a flat $n$-dimensional torus).  Then
\begin{itemize}
\item[(1)] $\;b^{j}_{(2)}(\widetilde{M},[\theta])=0$ for $j=0,\ldots,n$.

\item[(2)]
$\;\alpha_{j}(\widetilde{M},[\theta])=\left\{\begin{array}{lcl}
\frac{n}{2} & \mbox{ if } & [\theta]=0, \\[+5pt]
\infty & \mbox{ if } & [\theta]\neq 0.
\end{array}\right.$

\item[(3)] $\;\sigma(\Delta_{\theta,j})=[\|\theta\|^{2},\infty)$.

\item[(4)] $\;\lambda_{0,j}(\Delta_{\theta})=\|\theta\|^{2}$.
\end{itemize}
\end{proposition}

\subsection{Nil manifolds}
Let $G$ be a connected and simply-connected nilpotent Lie group and
$\Gamma$ be a torsion-free, discrete, cocompact subgroup.  Then the space
of differential forms on a nilmanifold $\,_{\Gamma}\backslash G$ is isomorphic
to the space of $\Gamma$-invariant differential forms on $G$, that is,
$$
\Omega^{*}(\,_{\Gamma}\backslash G)\cong\Omega^{*}(G)^{\Gamma}
$$
Let ${\mathcal J}$ denote the Lie algebra of $G$.
The Lie algebra cochains are
$$
C^{*}({\mathcal J})=\Omega^*(G)^{G}\hookrightarrow\Omega^{*}(G)^{\Gamma}.
$$
Nomizu \cite{Nom} proved that the inclusion above
induces an isomorphism in
cohomology, that is,
$$
H^{*}({\mathcal J})\cong H^{*}(\,_{\Gamma}\backslash G)
$$
In other words elements of $H^{*}(\,_{\Gamma}\backslash G)$ are represented
by nowhere zero, $G$-invariant differential forms on $G$.
In particular, if $\alpha\in H^{1}(\,_{\Gamma}\backslash G)$ and $\alpha$
is not trivial, then $\alpha$ is represented by a nowhere
zero closed $G$-invariant 1-form on $G$, which induces
a nowhere zero closed 1-form $\theta$ on  $\Gamma\backslash G$.
The asymptotic Morse inequalities imply that
$$
b^{j}_{(2)}(\widetilde{\Gamma\backslash G} ,s[\theta])=0\;\;\mbox{ for all
}\;\;j=0,1,\ldots,\dim
G,\;\mbox{ and }\;s\gg 0.
$$
In fact, one sees in the proof of these inequalities
 that a spectral gap near
zero develops for $\Delta_{s\theta,j}$, for all $s\gg 0$, since $\theta$ is
nowhere zero.  So
$$
\alpha_{j}(\widetilde{\Gamma\backslash G},s[\theta])=\infty\;\;\mbox{ for }\;\;
j=0,\ldots,\dim G,\;\mbox{ and }\;s\gg 0.
$$
Summarizing, one has

\begin{proposition}.  Let $\theta$ be a closed 1-form on a
compact Nil manifold
$M$ of dimension $n$.  Then
\begin{itemize}
\item[(1)] $\;b^{j}_{(2)}(\widetilde{M},s[\theta])=0$ for
$j=0,\ldots,n$, and for all
$s\gg 0$.

\item[(2)] $\;\alpha_{j}(\widetilde{M},s[\theta])=\infty$ for
$j=0,\ldots,n$, and
for all $s\gg 0$.

\item[(3)] $\;\lambda_{0,j}(\Delta_{s\theta})>0$ for $j=0,\ldots,n$, and for all
$s\gg 0$.
\end{itemize}
\end{proposition}

\subsection{Mapping cylinders}
Let $p:M\rightarrow S^{1}$ be a closed $n$-dimensional
manifold which is a fiber
bundle over the circle $S^{1}$.  Since there is a nowhere zero closed
1-form $\theta_{0}$ on $S^1$, we can consider $p^{*}(\theta_{0})$ which is
a nowhere zero closed 1-form on $M$.  Arguments similar to that given for
Nil manifolds yield

\begin{proposition} For
%Let $p^{*}(\theta_{0})$ be a nowhere zero closed 1-form on
a closed $n$-dimensional manifold $p:M\rightarrow S^{1}$ which is a
fiber bundle over $S^{1}$, in the notations described above one has
\begin{itemize}
\item[(1)] $\;b^{j}_{(2)}(\widetilde{M},s[p^{*}(\theta_{0})])=0$ for
$j=0,\ldots,n$,
and all $s\gg 0$.

\item[(2)] $\;\alpha_{j}(\widetilde{M},s[p^{*}(\theta_{0})])=\infty$ for
$j=0,\ldots,n$,
and all $s\gg 0$.

\item[(3)] $\;\lambda_{0,j}(\Delta_{s,p^{*}(\theta_{0})})>0$ for
$j=0,\ldots,n$, and all $s\gg 0$.
\end{itemize}
\end{proposition}

\subsection{2-manifolds}
\begin{proposition}  Let $M$ be a compact 2-dimensional manifold.
Then the
twisted $L^{2}$ Betti functions $b^{j}_{(2)}(\widetilde{M},[\theta])$ are
constant on $H^1(M,{\mathbb R} )$ for $j=0,1,2$.  More precisely,
$b^{1}_{(2)}(\widetilde{M}, [\theta])=-{\chi}(M)$ and
$b^{0}_{(2)}(\widetilde{M},
[\theta])=b^{2}_{(2)}(\widetilde{M}, [\theta])=0$
.\end{proposition}

\begin{proof}
If $\operatorname{genus}(M)=0$, then $H^1(M,{\mathbb R} )=0$ and the result is
clear.  If $\operatorname{genus}(M)>0$, then by Proposition 3.12  one has
$$
b^{0}_{(2)}(\widetilde{M},[\theta])-b^{1}_{(2)}(\widetilde{M},[\theta])+
b^{2}_{(2)}(\widetilde{M},[\theta])={\chi}(M).
$$
Since  $b^{0}_{(2)}(\widetilde{M},[\theta])=
b^{2}_{(2)}(\widetilde{M},[\theta]) = 0$
by Lemma 5.1, one deduces the proposition.
\end{proof}

An immediate corollary is

\begin{corollary}  Let $\theta$ be a Morse 1-form on a compact
2-dimensional
manifold $M$.  Then $m_{1}(\theta)\geq-{\chi}(M)$.\end{corollary}

\end{document}